\documentclass[lettersize,journal]{IEEEtran}

\ifCLASSINFOpdf

\else

\fi

\usepackage{cite}
\usepackage{amsmath}
\usepackage{amsthm}
\usepackage{xcolor}
\usepackage{amssymb}
\usepackage{comment}
\usepackage{graphicx}
\usepackage{algorithm}
\usepackage{algorithmic}
\usepackage{subfig}

\usepackage[acronym]{glossaries}
\usepackage{glossaries-prefix}

\newtheorem{theorem}{Theorem}
\newtheorem{definition}{Definition}

\newtheorem{standingassumption}{Standing Assumption}
\newtheorem{assumption}{Assumption}
\newtheorem{lemma}{Lemma}
\newtheorem{proposition}{Proposition}
\newtheorem{corollary}{Corollary}
\theoremstyle{remark}
\newtheorem{remark}{Remark}

\DeclareMathOperator{\col}{col}
\DeclareMathOperator*{\argmin}{arg\,min}

\newcommand{\GP}[1]{\textcolor{blue}{#1}}

\newcommand{\mc}{\mathcal}

\makeglossaries
\newacronym{GNEP}{GNEP}{generalized Nash equilibrium problem}
\newacronym{SGNEP}{SGNEP}{stochastic GNEP}
\newacronym{GNE}{GNE}{generalized Nash equilibrium}
\newacronym{SGNE}{SGNE}{stochastic generalized Nash equilibrium}
\newacronym{NE}{NE}{Nash equilibrium}
\newacronym{DRO}{DRO}{distributionally robust optimization}
\newacronym{MIP}{MIP}{mixed-integer program}
\newacronym{iid}{i.i.d.}{independent and identically distributed}
\newacronym{wrt}{w.r.t.}{with respect to}
\newacronym{KKT}{KKT}{Karush-Kuhn-Tucker}
\newacronym{CC}{CC}{chance constraint}
\newacronym{DR}{DR}{distributionally robust}
\newacronym{CVaR}{CVaR}{conditional value-at-risk}
\newacronym{DRNEP}{DRNEP}{distributionally robust Nash equilibrium problem}
\newacronym{DRNE}{DRNE}{distributionally robust Nash equilibrium}
\newacronym{vDRGNE}{vDRGNE}{variational DRGNE}
\newacronym{vSGNE}{vSGNE}{variational SGNE}
\newacronym{SNE}{SNE}{stochastic Nash equilibrium}
\newacronym{SNEP}{SNEP}{stochastic Nash equilibrium Problem}
\newacronym{VI}{VI}{variational inequality}
\newacronym{DRVI}{DRVI}{distributionally robust variational inequality}
\newacronym{DRGNE}{DRGNE}{distributionally robust generalized Nash equilibrium}
\newacronym{DRGNEP}{DRGNEP}{distributionally robust generalized Nash equilibrium problem}

\allowdisplaybreaks

\begin{document}

\title{Data-driven Wasserstein distributionally robust generalized Nash equilibrium problems with heterogeneous uncertainty}

\author{Georgios Pantazis,~\IEEEmembership{}
        Barbara Franci,~\IEEEmembership{}
    Sergio Grammatico~\IEEEmembership{}
\thanks{G. Pantazis is with the Dynamics and Control section at the Department of Mechanical Engineering, Eindhoven University of Technology, The Netherlands ({\tt\small G.Pantazis33@tue.nl}).}
\thanks{B. Franci is with the Department of Mathematical Sciences, Politecnico di Torino, Italy, supported by “Programma Rita
Levi Montalcini” of the Italian Ministry of University and Research (53\_RID26FRABAR) ({\tt\small barbara.franci@polito.it}). }
\thanks{S. Grammatico is at Delft Center for Systems and Control, Delft University of Technology, The Netherlands ({\tt\small s.grammatico@tudelft.nl}).}}

\markboth{Journal of \LaTeX\ Class Files,~Vol.~14, No.~8, August~2015}%
{Shell \MakeLowercase{\textit{et al.}}: Bare Demo of IEEEtran.cls for IEEE Journals}

\maketitle

\begin{abstract}
We study stochastic generalized Nash equilibrium (GNE) problems with heterogeneous, stochastic uncertainty in the cost functions. To address the distributional uncertainty, we let each agent construct an individual Wasserstein ambiguity set around an empirical distribution obtained from private data. We establish finite-sample robustness guarantees for any \gls{DR} GNE. Focusing on variational equilibria, we derive a \gls{DR} \gls{VI} reformulation and establish asymptotic consistency with the nominal stochastic game.  Under mild conditions, we derive an exact finite-dimensional generalized Nash equilibrium
reformulation. We then quantify how heterogeneous
worst-case distributions affect monotonicity, and derive conditions
for uniqueness and equilibrium sensitivity. To compute a variational \gls{DR} GNE, we propose an inexact projected method and establish
convergence under vanishing oracle errors, an asymptotic error bound under bounded oracle errors, and linear convergence under exact
oracle evaluations. 
\end{abstract}


\IEEEpeerreviewmaketitle

\section{Introduction}
A variety of applications in smart grids, communication, and social networks involve self-interested interacting decision makers \cite{Saad_2012}. Such systems are often analyzed through game theory, which models agents that optimize coupled individual cost functions subject to shared operational constraints.  In the presence of coupling constraints, the appropriate solution concept for the considered game is the so-called \gls{GNE}. Unlike in a standard
Nash equilibrium problem, where each agent optimizes over an individual local
feasible set, independent of the decisions of the other agents, in a
\gls{GNEP} the feasible set of each agent depends on the decisions of
the other agents.
In practice, however, uncertainty can affect both performance and safety of such systems. 
Existing work on stochastic games typically relies on specific probabilistic \cite{Kouvaritakis} or structural \cite{Aghassi2006,FukuSOCCP} assumptions on the uncertainty set. Since the true uncertainty may not follow the postulated model, recent works have also explored data-driven approaches 
to obtain robustness certificates against strategic changes \cite{Pantazis2023_apriori,Fabiani2020b,Fele2021} while
other results rely on approximation schemes to estimate the expected cost rather than learning the distribution \cite{Lei_2022,Franci_2021_merely,Franci_2021}.

When the underlying probability distribution is unknown and only limited data are available, computing a robust generalized Nash equilibrium can be challenging. 
A recent development in optimization is to consider the so-called \gls{DR} approach which allows one to handle distributional uncertainty using ambiguity sets constructed from finite data \cite{shapiro_distributionally_2017,Parys2020data,Zhen2021AUT,Boskos_2024}. In this regard, ambiguity sets based on optimal transport, and in particular Wasserstein ambiguity sets, have received considerable attention
\cite{mohajerin_esfahani_data-driven_2018,Fournier_2023}. 
Distributionally robust approaches have also been developed in optimal control, particularly in model predictive control (MPC), with emphasis on computational tractability, closed-loop guarantees, and robustness to distributional uncertainty \cite{mcallister2025drmpc,wu2023ambiguity,mcallister2024inherent,falconi2025lqg}. The difficulty is further amplified in generalized game settings with
coupling constraints and heterogeneous uncertainty. In our setting, each
agent constructs an individual Wasserstein ambiguity set from its own
private data and robustifies its own expected cost with respect to this
set. These individually robustified objectives are coupled through the
agents' decisions and the shared constraints. Moreover, the worst-case
distributions, associated with the agents' robust costs, can differ across
agents and may also depend on the collective decision profile. Consequently,
monotonicity properties of the underlying nominal game are not, in general,
automatically inherited by the \gls{DR} game. This
heterogeneous structure motivates both the finite-sample analysis of the ambiguity sets and the study of monotonicity preservation under heterogeneous ambiguity.

Applications of \gls{DRO} include electricity market participation \cite{nasiri2022tactical} and optimal power flow \cite{guo2018data}. In a
game-theoretic setting, \cite{Liu2018} develops \gls{DR}
Nash and Stackelberg equilibrium models for supply chain
applications. These works, however, do not consider Nash games
with shared constraints in which each agent constructs an
individual data-driven Wasserstein ambiguity set for their own expected cost.
Closer to \gls{GNEP}s,
\cite{fabiani2023distributionally,wen2025distributionally}
consider generalized Nash games with Wasserstein \gls{DR}
chance constraints, so that distributional ambiguity affects the shared
coupling constraints rather than the agents' expected cost functions.
In contrast, in the present work each agent constructs an individual
data-driven Wasserstein ambiguity set for its own expected cost, while the
shared coupling constraints are deterministic.
More recent works address other related but distinct formulations.
In particular, \cite{liu2025distributionally} studies \gls{DR} equilibrium seeking under unknown dynamics via Gaussian process regression, while \cite{liu2026finitehorizon} develops a finite-horizon distributionally robust equilibrium-seeking framework for unknown coupled dynamics. Moreover,
\cite{liu2026distributed} considers generalized \gls{DR}
equilibrium seeking for stochastic dynamical games with nonlinear coupled
constraints and unknown time-varying disturbances.
\gls{DR} hierarchical interactions are considered in
\cite{shen2026strategic} through a bilevel game-theoretic formulation. Recent works have also considered data-driven Wasserstein
\gls{DRNEP} with heterogeneous or private data.
In particular, \cite{wang2025lagrangian} develops a Lagrangian approach
to heterogeneous-data Wasserstein \gls{DRNE}
seeking. 
These works concern Nash games with decoupled strategy sets, whereas
the present paper addresses \gls{GNEP}s with shared convex
constraints and develops finite-sample, consistency, sensitivity, and
monotonicity results for the resulting \gls{DRGNEP}.
A preliminary version of this work \cite{PantazisFranci2023} introduced the
data-driven Wasserstein Nash equilibrium framework with heterogeneous
ambiguity. The present manuscript substantially
extends that preliminary framework from Nash games with Cartesian
strategy sets to \gls{GNEP}s with shared convex
constraints. The resulting theoretical and algorithmic developments are
summarized below.

In the presence of coupling constraints, we distinguish between general
\gls{DRGNE} and their subclass, the so-called \gls{vDRGNE}. Our existence result establishes
the existence of a \gls{vDRGNE}, and hence also the existence of a \gls{DRGNE},
while the out-of-sample robustness certificates hold for any
\gls{DRGNE}. In contrast, our \gls{VI}
analysis, pseudogradient sensitivity bounds, asymptotic consistency results,
as well as Algorithm 2, concern \gls{vDRGNE}s. The main difference is that under standard regularity conditions, a \gls{vDRGNE} can be
characterized by a common multiplier associated with the shared
constraints, whereas a general \gls{DRGNE} may involve different multipliers
for different agents. This common multiplier structure enables a
\gls{VI} characterization over the collective feasible set
and motivates our later focus on \gls{vDRGNE}s for  analytical and
algorithmic results that exploit their structure.

\subsection{Main Contributions}    
Differently from previous works, our paper focuses on data-driven Wasserstein \gls{DRGNEP}s with heterogeneous ambiguity affecting the cost function of each agent. 
In particular, our contributions are the following:
\begin{enumerate}
\item We introduce a class of \glspl{SGNEP} with coupling convex constraints, where each agent constructs an individual Wasserstein ambiguity set from private data with an individual radius. This framework allows the ambiguity sets to differ across agents due to
heterogeneous sample sizes, confidence levels, or radii chosen by the agents. Furthermore, we establish existence of a \gls{vDRGNE} under
convexity, compactness, and continuity assumptions, thereby
also guaranteeing nonemptiness of the \gls{DRGNE} set.
	\item We derive collective finite-sample guarantees showing that, at any \gls{DRGNE}, each
agent's true expected cost is bounded by its robust surrogate with high
confidence. Concerning \gls{vDRGNE}s, we obtain a
pseudogradient sensitivity bound and almost-sure asymptotic
consistency with the nominal variational equilibrium set.
Under strong monotonicity, we further derive an equilibrium
distance bound involving the deviation from the unknown true
distributions, together with a high confidence result,
based on the finite-sample Wasserstein radii.
	\item Under conditions ensuring attainment of the Wasserstein dual
problems, we derive exact finite-dimensional reformulations
of the DRGNEP. For a common ambiguity set and a shared uncertain cost
term, we provide conditions under which the associated
extended saddle mapping is monotone. Finally, for heterogeneous ambiguity sets, we quantify the effect of decision-dependent worst-case distributions on the  monotonicity of the mapping.
\end{enumerate}
To compute \gls{DRGNE}, we consider two algorithmic schemes: i) a primal-dual method for the finite-dimensional reformulation, inspired by~\cite{Franci_2021_merely,malitsky_golden_2020}, used as a heuristic in the general heterogeneous setting, and ii) an inexact vDRGNE seeking method with convergence
under vanishing oracle errors, an asymptotic error bound under
bounded oracle errors, and linear convergence under exact oracle
evaluations. We illustrate our framework on a peer-to-peer electricity market and a heterogeneous Nash-Cournot game with coupling constraints.

\section{Problem formulation}\label{sec_problem}
 \subsection{Preliminaries}

We denote by $\|\cdot\|_q$ the $q$-norm and by
$\|\cdot\|$ the Euclidean norm. We use the standard definitions
of projection, Lipschitz continuity, monotonicity, strong
monotonicity, and convex subdifferentials. The residual is
$\operatorname{res}_{C}(x)
:=
\|x-\operatorname{proj}_{C}[x-F(x)]\|$.
For $f:\mathbb R^d\times\mathbb R^p\to\mathbb R$, $f^\ast$
denotes its partial conjugate with respect to the second
argument, and $\iota_C$ denotes the indicator function of $C$. For a set $S$, $\operatorname{co}(S)$ denotes its convex hull.  Let $\mathcal{P}(\Xi)$ denote the set of all probability distributions
supported on~$\Xi$. The set of all probability distributions defined on $\Xi \subseteq \mathbb{R}^p$ with a bounded first-order moment is
\begin{equation}
\mathcal{M}(\Xi)
= \textstyle
\left\{
\mathbb{P}\in\mathcal{P}(\Xi)
\,\middle|\,
\int_{\Xi}\|\xi\|\,\mathbb{P}(d\xi)<\infty \nonumber
\right\}.
\end{equation}

 We can now define the Wasserstein metric to quantify the distance between two probability distributions \cite{villani_topics_2016}.

\begin{definition} \label{Wasserstein}
	The Wasserstein distance $d_W: \mathcal{M}(\Xi) \times \mathcal{M}(\Xi) \rightarrow \mathbb{R}_{\geq 0}$  between two distributions $\mathbb{Q}_1, \mathbb{Q}_2 \in \mathcal{M}(\Xi)$ is defined as 
	\begin{align}
		d_W(\mathbb{Q}_1, \mathbb{Q}_2):=&\inf_{\Pi \in \mathcal{J}(\xi_1 \sim \mathbb{Q}_1, \xi_2 \sim \mathbb{Q}_2)}\int_{\Xi^2} \|\xi_1-\xi_2 \| \Pi(d\xi_1, d\xi_2) , \nonumber 
	\end{align}
	where  $\mathcal{J}(\xi_1 \sim \mathbb{Q}_1, \xi_2 \sim \mathbb{Q}_2)$ represents the set of joint probability distributions of the random variables $\xi_1$ and $\xi_2$ with marginals  $\mathbb{Q}_1$, $\mathbb{Q}_2$, respectively.
\end{definition}
The Wasserstein metric can be viewed as the optimal transport plan to fit the probability distribution $\mathbb{Q}_1$ to $\mathbb{Q}_2$. 
An alternative dual interpretation can be derived by the so-called Kantorovich-Rubinstein theorem \cite{villani_topics_2016}:

\begin{theorem}\label{KR_theorem} 
Let $\mathcal{F}$ be the space of all 1-Lipschitz continuous functions on $\Xi$. Then, for two distributions $\mathbb{Q}_1$, $\mathbb{Q}_2\in \mathcal{M}(\Xi)$ 
		\begin{align}
		d_W(\mathbb{Q}_1, \mathbb{Q}_2):= \sup\limits_{f \in \mathcal{F}}\{\mathbb{E}_{{\mathbb{Q}}_1}[f(\xi)]-\mathbb{E}_{{\mathbb{Q}}_2}[f(\xi)] \}.\nonumber 
	\end{align} 
\end{theorem}


\subsection{Heterogeneous distributionally robust Generalized Nash Equilibrium Problems}
Consider a set of agents with index set $\mathcal{N}=\{1, \dots, N\}$. The decision vector of each agent $i \in \mathcal{N}$ is denoted by $x_i=\col((x^{(j)}_i)_{j=1}^{n}) \in X_i \subseteq \mathbb{R}^{n}$, where $x^{(j)}_i, j=1,\dots,n$, denotes an element. Let $x_{-i}=\col((x_j)_{j=1, j \neq i}^N) \in X_{-i} \subseteq \mathbb{R}^{(N-1)n}$ be the vector of the decision variables of all agents except for agent $i$ and let $x=\col((x_i)_{i =1}^N)\in X=\prod_{i\in\mc N} X_i$ be the collective decision vector. 
The agents are coupled through the collective feasible set
$\mathcal X
:=
\{x\in X\mid\psi(x)\leq0\},
$
where $\psi:\mathbb R^{nN}\to\mathbb R^{\tilde{n}}$
is componentwise convex and continuous.
Note that affine equality constraints can be included by representing
$Ax=b$ as $Ax-b\leq0$ and $b-Ax\leq0$.  We assume that
$\mathcal{X}$ is nonempty. For each $i\in\mathcal{N}$ and a given
$x_{-i}$, each agent's feasible set is given by
$\mathcal X_i(x_{-i})
:=\{
z_i\in X_i
\mid
(z_i,x_{-i})\in\mathcal X
\}$.

\begin{standingassumption}
For every $i\in\mathcal N$, the true distribution
$\mathbb P_i$ belongs to $\mathcal M(\Xi_i)$.
Moreover, for every $x\in\mathcal X$, the cost function
$h_i(x,\cdot)$ is measurable and $\mathbb P_i$-integrable,
so that $\mathbb E_{\mathbb P_i}[h_i(x,\xi_i)]$
is finite and well-defined. \hfill $\square$
\end{standingassumption}

Now we can consider the following game, where each agent \(i \in \mathcal{N}\) aims at solving 
\begin{align}
\underset{x_i\in\mathcal X_i(x_{-i})}{\operatorname{minimize}}
\quad
\mathbb{E}_{\mathbb{P}_i}
\bigl[h_i(x_i,x_{-i},\xi_i)\bigr],
\label{eq:nominal}
\end{align}
where 
	 $\mathbb{E}_{\mathbb{P}_i}[h_i(x_i, x_{-i}, \xi_i)]:= \int_{\Xi_i}h_i(x_i, x_{-i}, \xi_i) \mathbb{P}_i(d\xi_i).$
The parameter $\xi_i \in \Xi_i$ denotes the uncertainty affecting the cost of agent $i \in \mathcal{N}$. The  collection of stochastic optimization problems in \eqref{eq:nominal} constitutes a 
\gls{SGNEP} and its most common solution concept is known as a 
\gls{SGNE} \cite{Shanbhag_2011,Franci_2021_merely,Franci_2021}.

\begin{definition} \label{SNE}
A decision vector $x^*=(x^\ast_i)_{i \in \mathcal{N}} \in \mathcal{X}$ is a SGNE of the game in \eqref{eq:nominal} if
	\begin{equation}
		\forall  i \in \mathcal{N}: \quad	x_i^* \in \argmin_{x_i \in \mathcal{X}_i(x^\ast_{-i})} \mathbb{E}_{\mathbb{P}_i}[h_i(x_i, x^*_{-i}, \xi_i)].  \nonumber
	\end{equation}
\end{definition}
We assume that the probability distributions $\mathbb{P}_i$ of the uncertainties $\xi_i \in \Xi_i$ affecting each agent $i \in \mathcal{N}$ are unknown. 
Therefore, the formulation of \eqref{eq:nominal} and its corresponding solution concept have to be adapted to account for the possible ambiguity in the probability distributions. To this end,  we consider a set of possible distributions that each agent's uncertainty might follow, which contains the true probability distributions with high confidence (see Lemma \ref{guarantees}), the so-called \emph{ambiguity set}. To construct an individual ambiguity set for each agent  $i \in \mathcal{N}$, we follow a data-driven approach and calculate the empirical distribution $\mathbb{\hat{P}}_{K_i}$ from $K_i$ i.i.d. samples $\xi_{K_i}=\{\xi^{(k_i)}\}_{k_i=1}^{K_i} \in \Xi_i^{K_i}$ according to
\begin{align}
	\textstyle \mathbb{\hat{P}}_{K_i} = \frac{1}{K_i} \sum_{k_i=1}^{K_i} \delta(\xi- \xi^{(k_i)}). \nonumber 
\end{align}
Then, we construct the Wasserstein ambiguity set 
$$\mathbb{B}_{\epsilon_i}(\mathbb{\hat{P}}_{K_i}) := \{ \mathbb{Q}_i \in \mathcal{M}(\Xi_i): d_W(\mathbb{\hat{P}}_{K_i}, \mathbb{Q}_i) \leq \epsilon_i \}$$ with radius $\epsilon_i \in \mathbb{R}_{\geq 0}$ and centre the empirical probability distribution  $\mathbb{\hat{P}}_{K_i} \in \mathcal{M}(\Xi_i)$. 
\begin{remark}
The radius $\epsilon_i$ can have different interpretations. In the
statistical setting of Section~III, it can be calibrated as a
function of the sample size $K_i$ and confidence parameter $\beta_i$ so
that the ambiguity set contains the true distribution with a prescribed
high probability. Consequently, heterogeneous radii may arise from
different sample sizes or confidence requirements across agents.
Alternatively, an agent may choose a larger radius than the statistically
calibrated value to impose additional robustness against distributional
misspecification.  \hfill $\square$
\end{remark}
The problem in \eqref{eq:nominal} then becomes a data-driven \gls{DR} game:
\begin{equation} \label{eq:dr}
	\ \forall  i \in \mathcal{N}: \min_{x_i \in \mathcal{X}_i(x_{-i})} \sup_{\mathbb{Q}_i \in \mathbb{B}_{\epsilon_i}(\mathbb{\hat{P}}_{K_i})} \mathbb{E}_{\mathbb{Q}_i}[h_i(x_i, x_{-i}, \xi_i)].
\end{equation}

For each agent $i\in\mathcal N$, we define the robust cost function
\begin{equation}\label{eq:robust_objective}
\mathcal{R}_i(x;\hat{\mathbb{P}}_{K_i},\epsilon_i)
:=
\sup_{Q_i\in \mathbb{B}_{\epsilon_i}(\hat{ \mathbb{P}}_{K_i})}
\mathbb E_{Q_i}[h_i(x,\xi_i)].  
\end{equation}
For a fixed collection of datasets and ambiguity radii,
we write $\mathcal{R}_i(x)$ for brevity. Thus, each agent's problem
in (2) becomes
\begin{align}
\min_{x_i\in X_i(x_{-i})} \mathcal{R}_i(x_i,x_{-i}). \label{eq:R_i}
\end{align}
Although $h_i(\cdot,x_{-i},\xi_i)$ may be smooth,
$\mathcal{R}_i(\cdot,x_{-i})$ need not be differentiable, since it is
a supremum of expected cost functions. This motivates
the use of partial subdifferentials in the variational
formulation developed below.

Accordingly, the equilibrium concept needs to be extended to  the \emph{Wasserstein \gls{DRGNE}}:
\begin{definition} \label{DRNE}
For a multi-sample $\xi_{K_i}\in\Xi_i^{K_i}$ drawn by agent $i\in\mathcal N$, a point $x_K^\ast=\operatorname{col}\big((x_{i,K}^\ast)_{i=1}^N\big)\in \mathcal{X}$, where $K:=\sum_{i\in\mathcal N}K_i$,
is a Wasserstein DRGNE of the game in~\eqref{eq:dr} if for each $i\in\mathcal N$ it holds that:
\begin{align}
x_{i,K}^\ast\in
\arg\min_{x_i\in \mathcal{X}_i(x^\ast_{-i,K})}
\sup_{\mathbb{Q}_i\in \mathbb{B}_{\epsilon_i}(\hat{\mathbb{P}}_{K_i})}
\mathbb E_{\mathbb{Q}_i}\!\left[
h_i(x_i,x_{-i,K}^\ast,\xi_i)
\right]. \nonumber
\end{align}
\end{definition}

The subscript $K$ in $x_K^\ast$ and $x^\ast_{-i,K}$ denotes the dependence of the solution on the multi-sample $\xi_K$. Our formulation allows heterogeneous ambiguity sets and deterministic shared constraints. When $\mathcal{X}=X$, Definition 3 reduces to the data-driven DRNE considered in games without coupling constraints.

\subsection{Variational inequalities for \gls{DR} generalized games}

In this section, we develop a \gls{VI} framework to analyze \gls{DR} generalized games. This approach allows us to characterize the \gls{DRGNE} using tools from \gls{VI} theory, which provides both theoretical and computational advantages. We first impose the following assumptions. 
\begin{assumption}\label{ass:existence}
For every $i\in\mathcal N$, $X_i$ is nonempty, compact,
and convex. The mapping
$\psi:\mathbb R^{nN}\to\mathbb R^{\tilde n}$
is componentwise convex and continuous, and $\mathcal X:=\{x\in X\mid\psi(x)\leq0\}$
is nonempty.
\hfill $\square$
\end{assumption} 
Conditions for the attainment of worst-case distributions over
Wasserstein ambiguity sets are characterized in~\cite{gao_distributionally_2022}. 
\begin{assumption}\label{convex}
For every $i\in\mathcal N$ and every fixed
$x_{-i}\in X_{-i}$ and $\xi_i\in\Xi_i$,
$h_i(\cdot,x_{-i},\xi_i)$ is convex on an open convex
neighborhood of $X_i$. Moreover, there exists an open
neighborhood $U$ of $X$ such that each robust cost
$\mathcal R_i$ is finite and continuous on $U$.
\hfill $\square$
\end{assumption}
\begin{assumption}\label{lipschitz}
For each $i\in\mathcal N$ and $\xi_i\in\Xi_i$, the function
$h_i(\cdot,\xi_i)$ is continuously differentiable on the open neighborhood $U$ of $\mathcal X$ from Assumption \ref{convex}. Moreover, 
for each $i\in\mathcal N$, there exist constants $L_{i,j}\geq0$, with
$j=1,\ldots,n$, such that:
\[
\textstyle \left|
\frac{\partial h_i(x,\xi_i)}{\partial x_i^{(j)}}-
\frac{\partial h_i(x,\xi_i')}{\partial x_i^{(j)}}
\right|
\leq
L_{i,j}\|\xi_i-\xi_i'\|,
\]
for all $x\in U$, $\xi_i,\xi_i'\in\Xi_i$, and
$j=1,\ldots,n$. \hfill $\square$
\end{assumption}
\begin{assumption}\label{cont_diff}
For each $i\in\mathcal N$, the function
$h_i:\mathbb R^{nN}\times\Xi_i\to\mathbb R$ is measurable in
$\xi_i$ and continuously differentiable in $x_i$ for every
$(x_{-i},\xi_i)\in X_{-i}\times\Xi_i$. Moreover, for every $i\in\mathcal N$, the mapping $x\mapsto \mathbb E_{P_i}
[\nabla_{x_i}h_i(x,\xi_i)]$
is well defined and continuous on $\mathcal X$.\hfill $\square$
\end{assumption}
By Assumptions~\ref{lipschitz} and~\ref{cont_diff},
differentiation and expectation can be interchanged.
Indeed, for any fixed $\bar\xi_i\in\Xi_i$, Assumption~\ref{lipschitz}
implies, locally uniformly in $x\in U$,
\[ \textstyle 
\left|
\frac{\partial h_i(x,\xi_i)}{\partial x_i^{(j)}}
\right|
\leq
\left|
\frac{\partial h_i(x,\bar\xi_i)}{\partial x_i^{(j)}}
\right|
+
L_{i,j}\|\xi_i-\bar\xi_i\|,
\qquad j=1,\ldots,n.
\]
Note that, since $\mathbb Q_i\in\mathcal M(\Xi_i)$ has a finite
first-order moment, the right-hand side is locally
dominated by a $\mathbb Q_i$-integrable function.
If $h_i(x_i,x_{-i},\cdot)$ is also
$\mathbb Q_i$-integrable at  $x$,
differentiation with respect to $x_i$ can therefore
be interchanged with expectation locally, i.e. , $\nabla_{x_i}\mathbb E_{\mathbb Q_i}[h_i(x,\xi_i)]
=
\mathbb E_{\mathbb Q_i}[\nabla_{x_i}h_i(x,\xi_i)]$.
We define the nominal pseudogradient mapping as
\begin{align}
F_{\mathbb P}(x)
:=
\operatorname{col}\!\left(
\mathbb E_{P_i}[\nabla_{x_i}h_i(x,\xi_i)]
\right)_{i\in\mathcal N}. \nonumber 
\end{align} 

Under these assumptions, we focus on the variational equilibria of the
stochastic and \gls{DR} generalized Nash equilibrium
sets. In the presence of shared constraints, a generalized Nash
equilibrium is not, in general, equivalent to the solution of a
\gls{VI} over the collective feasible set. The
variational equilibrium constitutes a particular selection from the
\gls{GNE} set and, under the regularity conditions
specified below, is characterized by a common multiplier associated
with the shared constraints. A general GNE may instead admit multipliers specific to each agent.  To define the nominal variational equilibrium, consider the stochastic \gls{VI}, denoted by $\mathrm{VI}(\mathcal X,F_{\mathbb P})$: \begin{align}
\text{find } 
&x^\star\in\mathcal X
\nonumber\
\text{s.t} \ \  (y-x^\star)^\top F_{\mathbb P}(x^\star)\geq 0, \ \forall y\in\mathcal X,
\label{eq:nominal-vi}
\end{align}
Every solution of
$\mathrm{VI}(\mathcal X,F_{\mathbb P})$ is called a \gls{vSGNE}. When $\mathcal X=X$, the
variational equilibrium set coincides with the \gls{SNE} set.  The following lemma establishes a connection between the set of SNE and the solution set SOL($\mathcal{X}, F_\mathbb{P}$) of $\mathrm{VI}(X, F_\mathbb{P})$. A similar result is mentioned in \cite{Shanbhag_2011}. 
\begin{lemma}
\label{lem_original_VI}
Let Assumptions~\ref{ass:existence}-\ref{cont_diff} hold. Every point
$x^\star\in\operatorname{SOL}(\mathcal X,F_{\mathbb P})$ is a
vSGNE of the game in~\eqref{eq:nominal}.
If $\mathcal X=X$, then the vSGNE set coincides with the SNE set.
\end{lemma}
\begin{proof} The proof is obtained from
    \cite{Franci_2021_merely}. Indeed, for any
    $i\in\mathcal N$ and
    $z_i\in\mathcal X_i(x_{-i}^\star)$, the point
    $(z_i,x_{-i}^\star)$ belongs to $\mathcal X$. Therefore, the VI
    condition and convexity in $x_i$ imply that $x_i^\star$ is an
    optimal response over $\mathcal X_i(x_{-i}^\star)$. Hence,
    $x^\star$ is an SGNE. If $\mathcal X=X$, the converse follows from the
    Cartesian set characterization of Nash equilibria as solutions of
    the corresponding VI \cite[Proposition 1.4.2]{Pang1}.\end{proof}   

We now derive the corresponding VI formulation
for the \gls{DR} game. Under
Assumptions~\ref{ass:existence} and~\ref{convex}, the robust
cost $\mathcal R_i(\cdot,x_{-i})$ is convex.
A feasible profile $x^\star \in \mathcal{X}$ is a DRGNE if and only if,
for every $i\in\mathcal N$,
\[
0\in
\partial_{x_i}\mathcal R_i(x^\star)
+
N_{\mathcal{X}_i(x^\star_{-i})}(x_i^\star).
\]
Equivalently, for every $i\in\mathcal N$, there exists
$v_i^\star\in\partial_{x_i}\mathcal R_i(x^\star)$ such that $(v_i^\star)^\top(z_i-x_i^\star)\geq 0$,
 for all $z_i\in \mathcal{X}_i(x^\star_{-i})$.
Let us now define the set-valued pseudo-subgradient mapping as
\[
\mathcal F_K(x)
:=
\left\{
\operatorname{col}(v_1,\ldots,v_N)
\;\middle|\;
v_i\in\partial_{x_i}\mathcal R_i(x)
\text{ for every } i\in\mathcal N
\right\}.
\]
The set $\mathcal F_K(x)$ consists of all vectors obtained
by stacking partial subgradients of the agents' robust
costs. The \gls{DRVI},
denoted by $\operatorname{VI}(\mathcal X,\mathcal F_K)$,
is then defined as:
\[
\left\{
\begin{aligned}
&\text{find } x_K^\star\in\mathcal X
\text{ and } v_K^\star\in\mathcal F_K(x_K^\star)
\\
&\text{s.t.} \ \ \ 
(v_K^\star)^\top(y-x_K^\star)\geq 0, \ \forall y\in\mathcal X.
\end{aligned}
\right.
\]
As such, robust costs $\mathcal R_i$ need not be
differentiable, even when the underlying costs $h_i$
are smooth. Moreover, the suprema defining
$\mathcal R_i$ need not be attained. The set-valued
DRVI formulation accommodates both possibilities
without requiring the selection of worst-case
distributions. 

We denote the solution set of $\operatorname{VI}(\mathcal X,\mathcal F_K)$ by
\[
\operatorname{SOL}(\mathcal X,\mathcal F_K)
:=
\left\{
x\in\mathcal X
\;\middle|\;
\begin{array}{l}
\exists v\in\mathcal F_K(x)\text{ such that}\\
 v^\top(y-x)\geq 0,\
\forall y\in\mathcal X
\end{array}
\right\}.
\]
The points in $\operatorname{SOL}(\mathcal X,\mathcal F_K)$
are called \gls{vDRGNE}. Lemma~2 shows that vDRGNEs constitute a
subset of the DRGNE set.

\begin{lemma}
\label{DRNE-DRVI}
Let Assumptions~\ref{ass:existence} and~\ref{convex} hold.
Every point
$x_K^\ast\in\operatorname{SOL}(\mathcal X,\mathcal F_K)$
is a \gls{DRGNE} of~\eqref{eq:dr}. If $\mathcal X=X$, then
$\operatorname{SOL}(X,\mathcal F_K)$ coincides with the 
set of \gls{DRNE} of~\eqref{eq:dr}.
\hfill $\square$
\end{lemma}

\emph{Proof:}
See Appendix \ref{app_problem}. \hfill $\blacksquare$

We next characterize the robust partial subdifferentials
in terms of worst-case distributions. The following
additional assumption is used for this representation.
\begin{assumption}\label{assum:compact_uncertainty}
For any $i\in\mathcal N$, the uncertainty set
$\Xi_i\subseteq\mathbb R^p$ is compact and
$h_i$ is continuous on $U\times\Xi_i$, where $U$ is
the open neighborhood of $X$ in Assumption~\ref{convex}. 
\hfill $\square$
\end{assumption}
\begin{theorem}[Danskin's theorem {\cite[Theorem~4.13]{bonnans2000perturbation}}]
\label{thm:danskin}
Let $Z$ be compact, $\phi$ be continuous, and
$\phi(\cdot,z)$ be continuously differentiable with
continuous gradient. If
$f(x)=\max_{z\in Z}\phi(x,z)$ is convex, then
\[
\partial f(x)
=
\operatorname{co}
\{\nabla_x\phi(x,z)\mid z\in\arg\max_{z\in Z}\phi(x,z)\}.
\]
\end{theorem}
\begin{lemma}\label{lem:robust_subgradient}
Let Assumptions~2-4 and
\ref{assum:compact_uncertainty} hold. Then, for every
$i\in\mathcal N$ and $x\in \mathcal{X}$, it holds that:
\[
\partial_{x_i}\mathcal{R}_i(x)
=
\left\{
\mathbb E_{\mathbb Q_i}
[\nabla_{x_i}h_i(x,\xi_i)]
\;\middle|\;
\mathbb Q_i\in\mathcal A_{i,K}(x)
\right\},
\] where $\mathcal A_{i,K}(x)
:=
\arg\max_{\mathbb Q_i\in
\mathbb{B}_{\epsilon_i}(\hat{\mathbb P}_{K_i})}
\mathbb E_{\mathbb Q_i}[h_i(x,\xi_i)]$.
\end{lemma}
\emph{Proof:}
Fix $x_{-i}\in X_{-i}$ and define $\varphi_i(x_i,\mathbb Q_i)
:=
\mathbb E_{\mathbb Q_i}
[h_i(x_i,x_{-i},\xi_i)]$.
Since $\Xi_i$ is compact,
$\mathbb B_{\epsilon_i}(\hat{\mathbb P}_{K_i})$ is compact
with respect to $d_W$. By
Assumption~\ref{assum:compact_uncertainty},
$\varphi_i$ is continuous, and hence
$\mathcal A_{i,K}(x)$ is nonempty and compact.
Moreover, by Assumptions~2-4,
$\varphi_i(\cdot,\mathbb Q_i)$ is convex and continuously
differentiable, with $\nabla_{x_i}\varphi_i(x_i,\mathbb Q_i)
=
\mathbb E_{\mathbb Q_i}
[\nabla_{x_i}h_i(x,\xi_i)].$
Danskin's theorem therefore yields
\[
\partial_{x_i}\mathcal R_i(x)
=
\operatorname{co}
\left\{
\mathbb E_{\mathbb Q_i}
[\nabla_{x_i}h_i(x,\xi_i)]
\mid
\mathbb Q_i\in\mathcal A_{i,K}(x)
\right\}.
\]
Finally, $\mathcal A_{i,K}(x)$ is convex because
$\mathbb B_{\epsilon_i}(\hat{\mathbb P}_{K_i})$ is convex
and the objective is linear in $\mathbb Q_i$.
Since the map
$\mathbb Q_i\mapsto
\mathbb E_{\mathbb Q_i}[\nabla_{x_i}h_i(x,\xi_i)]$
is also linear, the set of active gradients is convex.
Hence the convex hull above is redundant. \hfill $\blacksquare$

Under the assumptions of
Lemma~\ref{lem:robust_subgradient}, the DRVI is therefore
equivalent to finding $x_K^\star\in\mathcal X$ and
$\mathbb Q_{i,K}^\star\in\mathcal A_{i,K}(x_K^\star)$ for each
$i\in\mathcal N$, such that:
\[
\sum_{i\in\mathcal N}
(y_i-x_{i,K}^\star)^\top
\mathbb E_{\mathbb Q_{i,K}^\star}
[\nabla_{x_i}h_i(x_K^\star,\xi_i)]
\geq0,
\ \forall y\in\mathcal X.
\]

If, in addition, $\psi$ is continuously differentiable
and the normal cone identity 
\begin{align}
N_{\mathcal X}(x)
=
N_X(x)+
\left\{
\nabla\psi(x)^\top\lambda
\;\middle|\;
\lambda\geq0,\ \lambda^\top\psi(x)=0
\right\}
\end{align}
holds for every $x\in\mathcal X$, a \gls{vDRGNE} is equivalently
characterized by a pseudo-subgradient
$v_K^\star=\operatorname{col}((v_{i,K}^\star)_{i\in\mathcal N})
\in\mathcal F_K(x_K^\star)$ and a common multiplier
$\lambda_c^\star\geq0$ such that, for every $i\in\mathcal N$,
\[
0\in v_{i,K}^\star
+N_{X_i}(x_{i,K}^\star)
+\nabla_{x_i}\psi(x_K^\star)^\top\lambda_c^\star,
\]
with $x_K^\star\in X$, $
\psi(x_K^\star)\leq0$, and
$(\lambda_c^\star)^\top\psi(x_K^\star)=0$.
A general \gls{DRGNE} may instead have different multipliers
for each agent associated with the shared constraints. 

We first show the following existence result.
\begin{proposition} \label{prop:existence}
Let Assumptions~\ref{ass:existence} and~\ref{convex} hold.
Then,  
$\operatorname{VI}(\mathcal X,\mathcal F_K)$ admits a solution. \hfill $\square$
\end{proposition}
\emph{Proof:}
Since every $X_i$ is compact, the Cartesian set
$X=\prod_{i\in\mathcal N}X_i$ is compact. Moreover, continuity of $\psi$ implies that $\mathcal X$
is closed. Hence, by
Assumptions~\ref{ass:existence} and~\ref{convex}, $\mathcal X$ is
nonempty, compact, and convex. For every $i\in\mathcal N$, the function
$\mathcal R_i(\cdot,x_{-i})$ is convex as the pointwise
supremum of convex expected value functions. Since
$\mathcal R_i$ is finite-valued and continuous on an open
neighborhood of $X$, its partial subdifferential
$\partial_{x_i}\mathcal R_i(x)$ is nonempty, compact, and convex
for every $x\in\mathcal X$. Joint continuity and convexity in
$x_i$ imply that
$x\mapsto\partial_{x_i}\mathcal R_i(x)$ has a closed graph and is
locally bounded. Therefore, it is upper semicontinuous on
$\mathcal X$, and $\mathcal F_K(x)
=
\prod_{i\in\mathcal N}
\partial_{x_i}\mathcal R_i(x)$
is upper semicontinuous with nonempty compact convex values. Since $\mathcal X$ is compact and $\mathcal F_K$ is locally
bounded with a closed graph, the set
$\mathcal V
:=
\operatorname{co}
\left(
\bigcup_{x\in\mathcal X}\mathcal F_K(x)
\right)$
is nonempty, compact, and convex. Consider now the mapping
$\Phi:\mathcal X\times\mathcal V
\rightrightarrows\mathcal X\times\mathcal V$:
\begin{align}
\Phi(x,v)
:=
\{
\operatorname{proj}_{\mathcal X}(x-v)
\}
\times\mathcal F_K(x). \nonumber 
\end{align}
The mapping $\Phi$ is upper semicontinuous with nonempty
compact convex values.
Thus, the Kakutani fixed-point theorem yields
$(x_K^\star,v_K^\star)\in\mathcal X\times\mathcal V$
such that $x_K^\star=
    \operatorname{proj}_{\mathcal X}
    (x_K^\star-v_K^\star)$ with $v_K^\star \in\mathcal F_K(x_K^\star).$
Since $\mathcal X$ is nonempty, closed, and convex, the
characterization of the Euclidean projection gives
\begin{align}
    \bigl(
    x_K^\star-v_K^\star-x_K^\star
    \bigr)^\top
    (y-x_K^\star)
    \leq 0, \forall y\in\mathcal X.
\end{align}
Since
$x_K^\star-v_K^\star-x_K^\star=-v_K^\star$,
this is equivalent to  $(y-x_K^\star)^\top v_K^\star
    \geq 0, \forall y\in\mathcal X$, i.e., 
$x_K^\star\in
\operatorname{SOL}(\mathcal X,\mathcal F_K)$. \hfill $\blacksquare$

Consequently, the game in~\eqref{eq:dr} admits at least one \gls{vDRGNE}. The variational characterization above provides a convenient representation of the original DRGNEP and allows us to
study its equilibria through the corresponding robust pseudogradient.
We next investigate the statistical reliability of DRGNEs constructed
from finitely many private samples and their relation to the equilibria of the underlying stochastic game.

\section{Statistical Guarantees and Consistency of Data-Driven DRGNEs}\label{sec_prop}

We first establish high-confidence out-of-sample guarantees on the true expected costs at every data-driven \gls{DRGNE}. For the set of vDRGNE, we then bound the deviation between the robust and nominal pseudogradients and use this bound to prove almost-sure consistency with the nominal vSGNE set.

\subsection{Out-of-Sample Certificates}\label{sec_oos}

In this section, we provide robustness certificates on the out-of-sample performance of the \gls{DRGNE}. Our main result is based on the works \cite{fournier_rate_2015,Fournier_2023}. 
For the sake of our analysis, let us then impose the following assumption:
\begin{assumption}\label{light_tail}
For each agent $i \in \mathcal{N}$, the true probability distribution
$\mathbb{P}_i$ belongs to $\mathcal{M}(\Xi_i)$. Furthermore, there exists
$a_i>1$ such that $A_i := \mathbb{E}_{\mathbb{P}_i}
\left[\exp\left(\|\xi_i\|^{a_i}\right)\right] < \infty$.
\hfill $\square$
\end{assumption}
Assumption \ref{light_tail} ensures that the probability distributions do not have heavy tails, which is often common in sample average approximation. 
We can now quantify the statistical reliability of the ambiguity set constructed from finite data.
\begin{lemma}
\label{guarantees}
	Let Assumption \ref{light_tail} hold. Then, the following confidence bound holds for each agent $i \in \mathcal{N}$:
	\begin{align}
		\mathbb{P}_i^{K_i}\{\xi_{K_i} \in \Xi^{K_i} : \mathbb{P}_i \in \mathbb{B}_{\epsilon_i}(\mathbb{\hat{P}}_{K_i})\} \geq 1-\beta_i,
	\end{align}
	where 
	\begin{align} \label{beta_i}
		\beta_i= \begin{cases} 
			c_i \exp(-b_iK_i\epsilon_i^{\max{\{p,2\}}}) & \text{ if } \epsilon_i \leq 1 \\ 
			c_i \exp(-b_iK_i\epsilon_i^{a_i})  &  \text{ if } \epsilon_i >1, 
		\end{cases}
	\end{align}
	$K_i \geq 1$, $p \neq 2$, and $\epsilon_i>0$; $c_i$ and $b_i$ are positive constants that only depend on  $a_i, A_i$ and $p$. \hfill $\square$
\end{lemma}
\begin{proof}
The result follows by an application of \cite[Theorem 2]{fournier_rate_2015} for each agent $i \in \mathcal{N}$. 
\end{proof}
The case $p=2$, which is critical for the $1$-Wasserstein distance, requires the logarithmic correction in~\cite[Theorem~2]{fournier_rate_2015}. The bound in~(\ref{beta_i}) connects the confidence $1-\beta_i$ to the radius $\epsilon_i$ and sample size $K_i$, and decays exponentially in $K_i$. For $\epsilon_i\leq1$, its dependence on $\epsilon_i$ weakens as the dimension $p$ increases (See~\cite{Fournier_2023} for details on the constants in~(\ref{beta_i}) and~(\ref{epsilon_i})).

As a consequence of Lemma \ref{guarantees}, we can derive a rule for tuning the radius of the Wasserstein balls. In particular, solving with respect to $\epsilon_i$ for $p\neq 2$, we obtain:
\begin{align} 
	\label{epsilon_i}
	\epsilon_i(K_i, \beta_i) = \begin{cases} 
		\left(\frac{\ln\left(\frac{c_i}{\beta_i}\right)}{b_iK_i}\right)^{\frac{1}{\max\{p,2\}}} & \text{if } K_i \geq \frac{\ln\left(\frac{c_i}{\beta_i}\right)}{b_i},  \\ 
		\left(\frac{\ln\left(\frac{c_i}{\beta_i}\right)}{b_iK_i}\right)^{\frac{1}{a_i}} & \text{otherwise}.  
	\end{cases} 
\end{align}
For the statistical guarantees developed in this
section, the Wasserstein radius is calibrated as
$\epsilon_i=\epsilon_i(K_i,\beta_i)$ according
to~\eqref{epsilon_i}. Evaluating this formula requires
valid concentration constants $b_i$ and $c_i$, whose
values depend on the tail parameters
of the considered distribution, as stated in Assumption \ref{light_tail}. Thus, the formula provides a theoretical
calibration, whose implementation, however,  requires
bounds on these parameters. Differences among
the agents' radii may then result from heterogeneous sample
sizes $K_i$ and/or confidence levels $\beta_i$.

Moreover, we provide a finite-sample certificate for the true
expected costs at any \gls{DRGNE}. Specifically, we show that with high confidence over
the samples, each agent's true expected cost is bounded
above by its corresponding \gls{DR} cost value, simultaneously for all agents. On the event that all ambiguity sets contain their respective
true distributions, these bounds hold for every feasible
decision profile and hence also for any \gls{DRGNE} computed
from the same samples. The following lemma formalizes this guarantee.

\begin{lemma}
\label{Lem:certificates_different} 
	Let $x^*_K$\footnote{Whenever a data-driven
equilibrium $x^*_K$ appears in a probability statement,
we assume that it is selected as a measurable
function of the datasets.} be a \gls{DRGNE} of the game in \eqref{eq:dr} and let $\mathcal{R}_i(x;\hat{\mathbb{P}}_{K_i},\epsilon_i)$ be as in \eqref{eq:robust_objective}.
	Then, for a drawn multisample $\xi_K=\{\xi_{K_i}\}_{i=1}^N$ it holds that	\begin{align}
			\mathbb{P}^K\{ \mathbb{E}_{\mathbb{P}_i}[h_i(x^*_K, \xi_i)] \leq \mathcal{R}_i(x^\ast_K;\hat{\mathbb{P}}_{K_i},\epsilon_i), \ \forall \  i  \in \mathcal{N} \} 	\nonumber \\ 
            \geq	1-\sum_{i \in \mathcal{N}}\beta_i, \nonumber 
		\end{align}  
	where $\mathbb P^K$ denotes the joint probability law of the
agents' datasets, with marginals $\mathbb P_i^{K_i}$
for each agent $i \in \mathcal{N}$\footnote{Note that each agent's $i \in \mathcal{N}$ dataset consists of i.i.d.\
samples from $\mathbb P_i$, while dependence across
agents' datasets is allowed for this result.}.
	\end{lemma}
	\begin{proof}
	See Appendix \ref{app_prop}. 
	\end{proof}
    If the agents' datasets are mutually independent,
the coverage events are independent, and the lower bound
in Lemma~\ref{Lem:certificates_different} improves to
$\prod_{i \in \mathcal{N}}(1-\beta_i)$, or to $(1-\beta)^N$ when
$\beta_i=\beta$ for all $i$. Finally, note that if all agents share the same true distribution,
multi-sample, and ambiguity radius, their coverage
events coincide. Thus, the collective expected cost bounds in
Lemma~\ref{Lem:certificates_different} hold with
confidence at least $1-\beta$.

\subsection{Sensitivity of the \gls{DR} Pseudogradient} \label{sec_asymp}

We next bound the deviation of the robust pseudogradient
$\mathcal F_K$ from its nominal counterpart $F_{\mathbb P}$.
This bound will subsequently be used to obtain asymptotic consistency of
the corresponding variational equilibrium sets.

\begin{definition}
	For a set $S\subseteq X$, let us define the deviation of the
set-valued mapping $\mathcal F_K$ from $F_\mathbb{P}$, parameterized by distribution $\mathbb{P}$, at $x\in S$ as $\operatorname{dev}_K(x)
:=
\sup_{v_K\in\mathcal F_K(x)}
\|v_K-F_\mathbb{P}(x)\|.$
We then define $\|\mathcal F_K-F_\mathbb{P}\|_S
:=
\sup_{x\in S}\operatorname{dev}_K(x).$
We say that $\mathcal F_K$ converges uniformly to $F_\mathbb{P}$ on $S$
if $\|\mathcal F_K-F_\mathbb{P}\|_S\to0 \text{ as }K_i\to\infty$, for all  $i\in\mathcal N.$  \hfill $\square$
\end{definition} 
The following result provides the bound between the \gls{DR} mapping and the nominal stochastic mapping. For the sequence of empirical games considered below, the neighborhood
in Assumption~\ref{convex} may depend on the dataset. The
constants $L_{i,j}$ in Assumption~\ref{lipschitz}, however, are properties of the cost functions and are therefore independent of $K$.

\begin{lemma}
\label{rho_upper_bound}
	Let Assumptions~1-4 hold and let us define $\tilde{w}_i=
\epsilon_i(K_i,\beta_i)
+d_W(\hat{\mathbb P}_{K_i},\mathbb P_i)$ and $\ell_{i,K}:=
\big(\sum\limits_{j=1}^{n}L_{i,j}^2\big)^{1/2}\tilde w_i$. Then, it holds that $ \|\mathcal F_K-F_{\mathbb P}\|_{\mathcal X} \leq\sqrt{\rho_K}$,
where $\rho_K:=
\sum\limits_{i=1}^N\sum\limits_{j=1}^n
L_{i,j}^2
\tilde{w}_i^2.$
 \hfill $\square$
\end{lemma}
\
\emph{Proof.}
For each $i\in\mathcal N$ and a given decision  $x\in\mathcal X$, we choose a sufficiently small $r>0$ such that
$(z_i,x_{-i})\in U$ whenever $\|z_i-x_i\|<r$,
and $\mathcal R_i(\cdot,x_{-i})$ is convex on $\mathbb{B}(x_i,r)$.
For any such $z_i$, let us now define $\Delta h_i(\xi_i)
:=h_i(x_i,x_{-i},\xi_i)-h_i(z_i,x_{-i},\xi_i)$, where the dependence on $x_i$, $z_i$, and $x_{-i}$
is suppressed for brevity.
By the fundamental theorem of calculus and Assumption~3 it then holds that:
\[
\begin{aligned}
&|\Delta h_i(\xi_i)-\Delta h_i(\xi_i')| \le
\left(\sum_{j=1}^{n}L_{i,j}^2\right)^{1/2}
\|x_i-z_i\|\|\xi_i-\xi_i'\|.
\end{aligned}
\]
Thus, by Theorem \ref{KR_theorem}, for every
$\mathbb Q_i\in \mathbb{B}_{\epsilon_i}(\hat{\mathbb P}_{K_i})$, it holds that:
\[
\begin{aligned}
\left|
\mathbb E_{\mathbb Q_i}[\Delta h_i(\xi_i)]
-\mathbb E_{\mathbb P_i}[\Delta h_i(\xi_i)]
\right| & \le
\left(\sum_{j=1}^{n}L_{i,j}^2\right)^{1/2}
\tilde w_i\|x_i-z_i\|\\
& =\ell_{i,K}\|x_i-z_i\|,
\end{aligned}
\]
where the last inequality follows from
\[
\begin{aligned}
d_W(\mathbb Q_i,\mathbb P_i)
&\le\epsilon_i(K_i,\beta_i)
+d_W(\hat{\mathbb P}_{K_i},\mathbb P_i) =\tilde w_i.
\end{aligned}
\]
For ease of notation, let us denote $\bar h_i(x):=\mathbb E_{\mathbb P_i}[h_i(x,\xi_i)]$
and $D_i(x):=\mathcal{R}_i(x)-\bar h_i(x)$.
For a fixed $\eta>0$, by the definition of the supremum, we consider a distribution
$\mathbb Q_i^\eta
\in \mathbb{B}_{\epsilon_i}(\hat{\mathbb P}_{K_i})$ such that $\mathbb E_{\mathbb Q_i^\eta}[h_i(x,\xi_i)]
\ge \mathcal{R}_i(x)-\eta$.
Then
\[
\begin{aligned}
D_i(x_i; x_{-i})-D_i(z_i;x_{-i}) &\le
\mathbb E_{\mathbb Q_i^\eta}[\Delta h_i(\xi_i)]
-\mathbb E_{\mathbb P_i}[\Delta h_i(\xi_i)]+\eta\\
&\le\ell_{i,K}\|x_i-z_i\|+\eta.
\end{aligned}
\]
Letting $\eta\downarrow0$ and interchanging $x_i$ and $z_i$ implies that $D_i(\cdot;x_{-i})$ is locally
$\ell_{i,K}$-Lipschitz. Now let $v_i\in\partial_{x_i}\mathcal{R}_i(x)$ and $\|d\|=1$.
For a sufficiently small $t>0$, convexity of
$\mathcal{R}_i(\cdot,x_{-i})$ and the Lipschitz bound above, it holds that:
\[
v_i^\top d
\le
\frac{\bar h_i(x_i+td,x_{-i})-\bar h_i(x_i,x_{-i})}{t}
+\ell_{i,K}.
\]
Letting $t\downarrow0$ and taking the supremum over
all directions $d$ with $\|d\|=1$, we obtain
$\|v_i-\nabla_{x_i}\bar h_i(x)\|\le\ell_{i,K}$
for all $i\in\mathcal N$.
Therefore, for every
$v_K=\operatorname{col}((v_i)_{i\in\mathcal N})
\in\mathcal F_K(x)$, it holds that:
\[
\|v_K-F_{\mathbb P}(x)\|^2
\le\sum_{i=1}^{N}\ell_{i,K}^2
=\rho_K.
\]
Taking the supremum over $v_K\in\mathcal F_K(x)$
and then over $x\in\mathcal X$ concludes the proof.
\hfill $\blacksquare$

\subsection{Asymptotic Consistency of vDRGNEs}
Under a suitable class of sequences of confidence parameters and Wasserstein radii, we show that the uniform deviation of the set-valued robust pseudogradient $\mathcal F_K$ from the nominal
pseudogradient $F_{\mathbb P}$ converges to zero almost surely. This result will then be used to establish the asymptotic
consistency of the vDRGNE set. For the asymptotic analysis of this section, let
$(\xi_i^{(k)})_{k\geq 1}$ be an infinite sequence of
i.i.d.\ samples from $\mathbb P_i$ for each agent 
$i\in\mathcal N$. The multi-sample of size $K_i$ consists
of the first $K_i$ samples of this sequence.
Let $\mathbb P^\infty$ denote the joint probability
distribution of these sequences. The induced joint distribution of the
corresponding finite datasets is denoted by
$\mathbb P^K$ and has marginals $\mathbb P_i^{K_i}$,
where $\mathbb P_i^{K_i}=\prod_{k=1}^{K_i} \mathbb{P}_i$. If the agents' sequences are mutually
independent, then
$\mathbb P^K=\prod_{i\in\mathcal N}\mathbb P_i^{K_i}$. We stress that
all almost-sure convergence statements below are
understood with respect to the probability measure $\mathbb P^\infty$.

\begin{lemma}\label{rho_convergence}
Let Assumption~\ref{light_tail} hold. For each $i\in\mathcal N$, let $(\beta_i^{(K_i)})_{K_i\in\mathbb N}\subset(0,1)$ be such that $\epsilon_i(K_i,\beta_i^{(K_i)})\rightarrow0 \text{ as }K_i\to\infty$. Then, $\rho_K\rightarrow0,  
\mathbb P^\infty\text{-almost surely }$,
as $\min_{i\in\mathcal N}K_i\to\infty$. \hfill $\square$
\end{lemma} 
\emph{Proof:} See Appendix \ref{app_prop}. 

We now establish that sequences of vDRGNEs approach
the nominal vSGNE set almost surely.
\begin{theorem}\label{asymptotic_consistency} Let Assumptions~1-4 and~6 hold  and let the Wasserstein radii satisfy the conditions of Lemma~\ref{rho_convergence}. For any sequence $x_K^\ast\in\operatorname{SOL}(\mathcal{X},\mathcal F_K)$, every accumulation point $\bar x$ of $(x_K^\ast)_K$ satisfies $\bar x\in\operatorname{SOL}(\mathcal{X},F_{\mathbb P})$ $\mathbb P^\infty$-almost surely, as $\min_{i\in\mathcal N}K_i\to\infty$. \hfill $\square$ \end{theorem}
\begin{proof}

Let $(x_K^\ast)_K$ be any sequence such that $x_K^\ast\in\operatorname{SOL}(\mathcal{X},\mathcal F_K)$. By the definition
of $\operatorname{SOL}(\mathcal{X},\mathcal F_K)$, for every $K$ there exists
$v_K^\ast\in\mathcal F_K(x_K^\ast)$ such that $(y-x_K^\ast)^\top v_K^\ast\geq0, \forall y\in \mathcal{X}$.
By Lemmas~\ref{rho_upper_bound} and~\ref{rho_convergence}, it holds that:
\begin{align}
\|v_K^\star-F_{\mathbb P}(x_K^\star)\|
&\leq
\|\mathcal F_K-F_{\mathbb P}\|_{\mathcal X}
\leq\sqrt{\rho_K}
\rightarrow 0,
\end{align}
$\mathbb P^\infty$-almost surely, as
$\min_{i\in\mathcal N}K_i\to\infty$. Since $\mathcal X$ is compact, the sequence
$(x_K^\ast)_K\subseteq\mathcal X$ admits an accumulation point.
Let $\bar x$ be an arbitrary accumulation point and select a
subsequence $(x_{K_r}^\ast)_r$ converging to $\bar x$.
By the continuity of $F_{\mathbb P}$ on $\mathcal{X}$,
$F_{\mathbb P}(x_{K_r}^\ast)
\rightarrow
F_{\mathbb P}(\bar x).$ From the triangle inequality it holds that:
\[
\begin{aligned}
\|v_{K_r}^\ast-F_{\mathbb P}(\bar x)\|
&\leq
\|v_{K_r}^\ast-F_{\mathbb P}(x_{K_r}^\ast)\| \\
&\quad+
\|F_{\mathbb P}(x_{K_r}^\ast)-F_{\mathbb P}(\bar x)\|
\rightarrow0.
\end{aligned}
\]
Hence,
$v_{K_r}^\ast\to F_{\mathbb P}(\bar x)$. Passing the limit in the relation $(y-x_{K_r}^\ast)^\top v_{K_r}^\ast\geq0, \forall y\in \mathcal{X}$, yields $(y-\bar x)^\top F_{\mathbb P}(\bar x)\geq0, \forall y\in \mathcal{X}$.
Therefore, $\bar x\in\operatorname{SOL}(\mathcal{X},F_{\mathbb P})$,
which concludes the proof.\end{proof} 

\begin{corollary}
Under the conditions of Theorem~\ref{asymptotic_consistency}, for any sequence $x_K^\ast\in\operatorname{SOL}(\mathcal{X},\mathcal F_K)$, $\operatorname{dist} \bigl(x_K^\ast,\operatorname{SOL}(\mathcal{X},F_{\mathbb P})\bigr) \rightarrow 0$, $\mathbb P^\infty$-almost surely as $\min_{i\in\mathcal N}K_i\to\infty$. If $\operatorname{SOL}(\mathcal{X},F_{\mathbb P})=\{x^\ast\}$, then $x_K^\ast\to x^\ast$ almost surely.
\end{corollary}

Theorem \ref{asymptotic_consistency} and Corollary 1 establish almost-sure
consistency in the sense that every accumulation
point of a sequence of vDRGNEs is a nominal vSGNE,
and the distance of the sequence to the nominal
vSGNE set converges to zero as
$\min_{i\in\mathcal N}K_i\rightarrow \infty$.
In case the nominal vSGNE is unique, the sequence converges
to that equilibrium almost surely.  Having established the statistical properties and asymptotic consistency of the \gls{DRGNEP}, we next derive finite-dimensional reformulations
that eliminate optimization over probability
distributions. Under additional structural
assumptions, the resulting supremum expressions then
admit explicit representations.

\section{Exact Finite-Dimensional Reformulations}\label{sec_finite}

\subsection{Wasserstein Dual Reformulation}
We now show that, under mild assumptions, the \gls{DRGNEP} in \eqref{eq:dr} is equivalent to a deterministic \gls{GNEP} with individual coupling constraints. The dualization below affects only the \gls{DR} part of the objective function, while the coupling constraint \(x_i\in\mathcal X_i(x_{-i})\) is not affected.

\begin{assumption} \label{pointwise_max}
For every $i\in\mathcal N$, the uncertainty set
$\Xi_i\subseteq\mathbb R^p$ is convex and closed. Moreover,
for every $(x,\xi_i)\in \mathcal{X}\times\Xi_i$, $h_i(x,\xi_i)
=
\max_{\ell_i=1,\ldots,L_i}
h_{\ell_i}(x,\xi_i)$,
where $-h_{\ell_i}(x,\cdot)$ is convex, proper, and lower
semicontinuous on $\Xi_i$, and $h_{\ell_i}(x,\cdot)$ is not
identically $-\infty$, for every
$\ell_i\in\{1,\ldots,L_i\}$.
\end{assumption}
Representing the cost function of each agent as a combination of pointwise maxima of affine functions with respect to the uncertainty can approximate a vast class of functions. 
Therefore, the following results are quite general.
\begin{theorem}
\label{reformulation_1}
	Let Assumption \ref{pointwise_max} hold. 
    For each agent $i \in \mathcal{N}$ and every $x \in \mathcal{X}$, the robust cost
$\mathcal R_i(x;\hat{\mathbb P}_{K_i},\epsilon_i)$
defined in \eqref{eq:robust_objective} admits the equivalent formulation:
\[
\mathcal R_i(x;\hat{\mathbb P}_{K_i},\epsilon_i)
=
\inf_{\lambda_i\geq0}
\left\{
\lambda_i\epsilon_i+
\frac{1}{K_i}\sum_{k_i=1}^{K_i}
g_i(x_i,x_{-i},\lambda_i,\xi_i^{(k_i)})
\right\}.
\]
If the infimum is attained for every  $x \in \mathcal{X}$, then the \gls{DRGNEP} in \eqref{eq:dr} admits the following reformulation: 
	\begin{align} 
		\forall	\  i \in \mathcal{N}: \min\limits_{x_i \in \mathcal{X}_i(x_{-i}), \lambda_i \geq 0} \lambda_i\epsilon_i + \frac{1}{K_i}\sum\limits_{k_i=1}^{K_i}  g_i(x_i,x_{-i},\lambda_i,\xi_i^{(k_i)}), \nonumber 
	\end{align}
	where
\begin{equation}
g_i(x_i,x_{-i},\lambda_i,\xi_i^{(k_i)})
=\sup\limits_{\xi_i\in\Xi_i}
(h_i(x_i,x_{-i},\xi_i)
-\lambda_i\|\xi_i-\xi_i^{(k_i)}\|). \nonumber 
\end{equation}
\end{theorem}
\begin{proof} See Appendix \ref{app_finite}.
\end{proof} \par 
The decision variable $\lambda_i \in \mathbb{R}_{\geq 0}$ represents the dual variable  of the Wasserstein ball constraint imposed on each probability distribution for each $i \in \mathcal{N}$. 
\subsection{Epigraphic GNEP Reformulation}
The following lemma then holds:
\begin{lemma} \label{reformulation_2}
		Let Assumption \ref{pointwise_max} hold. Suppose that,
for every $i\in\mathcal N$ and $x\in\mathcal X$,
the infimum over $\lambda_i$ in Theorem \ref{reformulation_1} is attained.
Then, the game in \eqref{eq:dr} admits the following equivalent reformulation, where each agent
$i\in\mathcal N$ solves:
	\begin{align} \label{eq:reformulation_convex}
    &\min_{x_i \in \mathcal{X}_i(x_{-i}), \lambda_i \geq 0, s_{k_i}, z_{k_i, \ell_i}} \lambda_i \epsilon_i + \textstyle \frac{1}{K_i}  \sum_{k_i=1}^{K_i} s_{k_i} & \nonumber \\
    & \  \mathrm{ s.t. }  \sup_{\xi_i  \in \Xi_i}[h_{\ell_i}(x_i, x_{-i}, \xi_i)-z_{k_i, \ell_i}^\top \xi_i] + z_{k_i, \ell_i}^\top \xi^{(k_i)} \leq s_{k_i}, \forall k_i, \forall \ell_i & \nonumber \\
    &\ \ \  \ \  \ \ \ \ \ \  \|z_{k_i, \ell_i}\| \leq \lambda_i,  \forall k_i, \forall \ell_i, \nonumber  
    \end{align}
where $k_i \in \{1, \dots, K_i\}$ and $\ell_i \in \{1, \dots, L_i\}$. 
\end{lemma}
\begin{proof} See Appendix \ref{app_finite}.
\end{proof}  

In Lemma \ref{reformulation_2}, $s_{k_i}$ denotes the epigraphic variable associated with a particular sample $k_i=1, \dots, K_i$, while $z_{k_i\ell_i}$ are auxiliary variables indexed by $\ell_i=1, \dots, L_i$, related to the elements of the pointwise maximum as per Assumption \ref{pointwise_max}.

For compactness of presentation, we denote $Y_i= X_i\times\mathbb R_{\geq 0} \times\mathbb R^{pK_iL_i}
\times\mathbb R^{K_i}$, $z_i=\col(((z_{k_i, \ell_i})_{k_i=1}^{K_i})_{\ell_i=1}^{L_i})$, $s_i=\col((s_{k_i})_{k_i=1}^{K_i})$, $y_i=\col(x_i, \lambda_i, z_i, s_i)$, $y_{-i}=\col(x_{-i}, \lambda_{-i}, z_{-i}, s_{-i})$ and $\phi_{k_i,\ell_i}(y_i,x_{-i})
:=
\bigl[-h_{\ell_i}(x_i,x_{-i},\cdot)
+\iota_{\Xi_i}\bigr]^*(-z_{k_i,\ell_i})
+
z_{k_i,\ell_i}^{\top}\xi_i^{(k_i)}
-s_{k_i}.$ Here, the conjugate function is taken with respect to the
uncertainty variable $\xi_i$, with $(x_i,x_{-i})$ fixed. Let
\[
\Phi_i(x_{-i})
:=
\bigcap_{\ell_i=1}^{L_i}
\bigcap_{k_i=1}^{K_i}
\left\{
y_i :
\phi_{k_i,\ell_i}(y_i,x_{-i})\leq 0,\;
\|z_{k_i,\ell_i}\|\leq\lambda_i
\right\},
\]
and define $H_i^{\mathcal{X}}(x_{-i}) := \left\{ y_i\in Y_i\cap\Phi_i(x_{-i})\mid (x_i,x_{-i})\in\mathcal{X} \right\}$.
According to Lemma \ref{reformulation_2}, the game in \eqref{eq:dr} takes the form:

\begin{equation} \label{eq:dr_eq}
\forall i\in\mathcal N:
\left\{
\begin{aligned}
\min_{y_i}\quad&
\lambda_i\epsilon_i+ \textstyle
\frac{1}{K_i}\sum_{k_i=1}^{K_i}s_{k_i}\\
\text{s.t.}\quad&
y_i\in H_i^{\mathcal{X}}(x_{-i}).
\end{aligned}
\right.
\end{equation} 
 The reformulation in ~\eqref{eq:dr_eq} replaces optimization
over probability distributions with a finite number of
decision variables. Its constraints nevertheless
retain suprema over the uncertainty support,
represented through convex conjugate functions. The polytopic support and
piecewise affine structure of the cost functions considered in the subsequent developments yield
an explicit tractable formulation of these constraints.
\begin{definition} \label{DRGNE}
A point $y_K^\ast
=
\operatorname{col}\bigl((y_i^\ast)_{i\in\mathcal N}\bigr)
=
(x_K^\ast,\lambda_K^\ast,z_K^\ast,s_K^\ast)$
is a GNE of the game in \eqref{eq:dr_eq} if, for every
$i\in\mathcal N$,
\[
y_i^\ast
\in
\arg\min_{y_i\in H_i^{\mathcal{X}}(x_{-i,K}^\ast)}
\left\{
\lambda_i\epsilon_i+ \textstyle
\frac{1}{K_i}\sum_{k_i=1}^{K_i}s_{k_i}
\right\}.
\]
\end{definition}
The following lemma connects the set of GNE of the game in \eqref{eq:dr_eq} with the \gls{DRGNE} of \eqref{eq:dr}.
\begin{lemma} 
\label{equilibrium_equivalence} 

Consider the setting of Lemma~\ref{reformulation_2}.  Then, if
 $y_K^\ast = (x_K^\ast,\lambda_K^\ast,z_K^\ast,s_K^\ast)$ is a GNE of the game in \eqref{eq:dr_eq}, then $x_K^\ast$ is a DRGNE of the game in \eqref{eq:dr}. 
\end{lemma}
\begin{proof}
By Lemma~\ref{reformulation_2}, each agent's problem in \eqref{eq:dr_eq} is an exact reformulation of its problem in \eqref{eq:dr}, with the deterministic feasible domain $\mathcal{X}_i(x_{-i})$ remaining unchanged. As such, the $x_i$-component of $y_i^\ast$ is a best response over $\mathcal{X}_i(x_{-i,K}^\ast)$, for every $i\in\mathcal N$. Therefore, $x_K^\ast$ satisfies Definition~3 and is a DRGNE of \eqref{eq:dr}. Conversely, let $x_K^\star$ be a DRGNE of the game in
\eqref{eq:dr}. Since the infimum is attained
(see Lemma~\ref{reformulation_2}), for each agent $i \in \mathcal{N}$
there exist feasible auxiliary variables 
$(\lambda_i^\star,z_i^\star,s_i^\star)$ such that $\lambda_i^\star\epsilon_i+
\frac{1}{K_i}\sum_{k_i=1}^{K_i}s_{k_i}^\star
=
\mathcal R_i(x_K^\star)$.
For any feasible unilateral deviation
$(x_i,\lambda_i,z_i,s_i)$ in the lifted game,
\[
\begin{aligned}
\lambda_i\epsilon_i+ \textstyle 
\frac{1}{K_i}\sum_{k_i=1}^{K_i}s_{k_i}
&\geq \mathcal R_i(x_i,x_{-i,K}^\star)\\
&\geq \mathcal R_i(x_K^\star)=\lambda_i^\star\epsilon_i+ \textstyle 
\frac{1}{K_i}\sum_{k_i=1}^{K_i}s_{k_i}^\star.
\end{aligned}
\]
The first inequality follows from the epigraphic reformulation, and the second from the definition of a  DRGNE.
Since the auxiliary variables are private to each
agent, their selection does not alter the other
agents' feasible deviations. Thus, the selected
lifted profile is a GNE.
\end{proof} \par  
\par 
\begin{remark}
Note that the converse establishes only the existence of auxiliary
variables, thus extending each DRGNE to an equilibrium
of the lifted game. These variables, however, are not specified
by the original equilibrium $x_K^\star$ alone and need not be unique. \hfill $\square$
\end{remark}
\subsection{Explicit Reformulation under Polytopic Support}
The following proposition focuses on the case where the support set $\Xi_i$ is polytopic and allows to obtain an explicit form for the conjugate function in each agent's constraints. 
\begin{proposition}\label{lem:polytopic}
Let each agent's support set $\Xi_i$ be a polytope, i.e.,
\[
\Xi_i
=
\left\{
\xi_i\in\mathbb{R}^{p}
\mid
C_i\xi_i\leq d_i
\right\},
\]
where $C_i\in\mathbb{R}^{m\times p}$ and
$d_i\in\mathbb{R}^{m}$. Moreover, let
\[
h_{\ell_i}(x,\xi_i)
=
a_{\ell_i}(x)^\top\xi_i+b_{\ell_i}(x),
\]
for every $\ell_i\in\{1,\ldots,L_i\}$ and
$i\in\mathcal N$. Then, a \gls{DRGNE} of \eqref{eq:dr}
can be obtained by solving the following \gls{GNEP}, for every
$i\in\mathcal N$:
\begin{equation}\label{eq_prop}
\begin{cases}
\displaystyle
\min_{\substack{x_i\in \mathcal{X}_i(x_{-i}),\;\lambda_i\geq0,\\
s_{k_i},\;\gamma_{k_i,\ell_i}}}
\quad
\lambda_i\epsilon_i
+ \textstyle 
\frac{1}{K_i}\sum_{k_i=1}^{K_i}s_{k_i}
\\[1ex]
\mathrm{s.t.}\quad
b_{\ell_i}(x_i,x_{-i})
+
a_{\ell_i}(x_i,x_{-i})^\top\xi_i^{(k_i)}
\\
\qquad\qquad
+
\gamma_{k_i,\ell_i}^\top
\left(d_i-C_i\xi_i^{(k_i)}\right)
\leq s_{k_i},
\quad
\forall k_i,\;\forall\ell_i,
\\
\qquad
\left\|
C_i^\top\gamma_{k_i,\ell_i}
-
a_{\ell_i}(x_i,x_{-i})
\right\|
\leq\lambda_i,
\quad
\forall k_i,\;\forall\ell_i,
\\
\qquad
\gamma_{k_i,\ell_i}\geq0,
\quad
\forall k_i,\;\forall\ell_i,
\end{cases}
\end{equation}
where $k_i\in\{1,\ldots,K_i\}$ and
$\ell_i\in\{1,\ldots,L_i\}$.
\end{proposition}
\begin{proof} See Appendix \ref{app_finite}.  \end{proof} 

\subsection{Reformulation under a Homogeneous Ambiguity Set}
We consider a special case of
Proposition~\ref{lem:polytopic} in which the agents share
both the ambiguity set and the uncertain part of their costs.
\begin{assumption}\label{ass:common_ambiguity}
All agents share the Wasserstein ambiguity set
$\mathbb B_{\epsilon}(\hat{\mathbb P}_M)$, where
\[
\hat{\mathbb P}_M
=
\textstyle  \frac{1}{M}\sum_{k=1}^{M}\delta_{\hat\xi^{(k)}}
\]
and the support $\Xi=\{\xi\in\mathbb R^p\mid C\xi\leq d\}$
is a nonempty polytope, with
$C\in\mathbb R^{m\times p}$ and $d\in\mathbb R^m$.
For every $i\in\mathcal N$, the cost has the form $h_i(x,\xi)=f_i(x)+g(x,\xi)$,
where $f_i$ is continuously differentiable and
\[
g(x,\xi)
=
\max_{\ell=1,\ldots,L}
\{a_\ell(x)^\top\xi+b_\ell(x)\}.
\]
Each $a_\ell$ is affine and each $b_\ell$ is finite
and convex on an open convex neighborhood of $X$.
\hfill $\square$
\end{assumption}
Under Assumption~\ref{ass:common_ambiguity}, the robust cost of each agent
$i\in\mathcal N$ is $\mathcal R_i(x)=f_i(x)+r_M(x)$, where
$r_M(x):=\sup_{\mathbb Q\in\mathbb B_\epsilon(\hat{\mathbb P}_M)}
\mathbb E_{\mathbb Q}[g(x,\xi)]$.
Thus, all agents have the same robust uncertain term, while their deterministic costs may differ.
Introduce the common auxiliary variables
$u=\operatorname{col}
\bigl(\lambda,(\gamma_{k\ell})_{k,\ell},s\bigr)$,
$s=\operatorname{col}(s_1,\ldots,s_M)$,
which take values in $\mathcal U
=
\mathbb R_{\geq0}
\times(\mathbb R_{\geq0}^{m})^{ML}
\times\mathbb R^M$. We can then define the cost
$J_M(u):=\lambda\epsilon+\frac{1}{M}\sum_{k=1}^{M}s_k$.
For $k=1,\ldots,M$ and $\ell=1,\ldots,L$, let us denote
\[
\begin{aligned}
c_{k\ell}^{(1)}(x,u)
&:=
b_\ell(x)+a_\ell(x)^\top\hat\xi^{(k)}+
\gamma_{k\ell}^\top(d-C\hat\xi^{(k)})-s_k,
\\
c_{k\ell}^{(2)}(x,u)
&:=
\|C^\top\gamma_{k\ell}-a_\ell(x)\|-\lambda.
\end{aligned}
\]
Finally, let $c_M(x,u)\in\mathbb R^{2ML}$ collect and stack these constraints for all values of $k$ and $\ell$.
By Proposition~\ref{lem:polytopic}, the common robust
term admits the value representation
\begin{equation}
r_M(x)
=
\min_{u\in\mathcal U}
\{J_M(u)\mid c_M(x,u)\leq0\}, \nonumber 
\end{equation}
while the physical shared constraints remain encoded in
$x\in\mathcal X$.
As such, the deterministic pseudogradient is given by $F_f(x)
:=
\operatorname{col}
\bigl((\nabla_{x_i}f_i(x))_{i\in\mathcal N}\bigr)$
while, for $\mu\in\mathbb R_{\geq0}^{2ML}$, the Lagrangian function takes the form:
\[
L_M(x,u,\mu)
:=
J_M(u)+\mu^\top c_M(x,u).
\]
For fixed $\mu$, let $\partial_{(x,u)}L_M$ denote
the joint subdifferential with respect to $(x,u)$.
The associated extended mapping is
\begin{equation}
\mathcal T_M(x,u,\mu)
:=
\begin{bmatrix}
\operatorname{col}(F_f(x),0)
+\partial_{(x,u)}L_M(x,u,\mu)
\\[1mm]
-c_M(x,u)
\end{bmatrix}. \nonumber 
\end{equation}
Here, $\partial_{(x,u)}L_M$ denotes the joint
subdifferential with respect to $(x,u)$, for fixed  $\mu$. As such, we note that the first block is in fact set-valued. The following result establishes the monotonicity of the augmented operator.
\begin{proposition}\label{prop:common_monotonicity}
Let Assumption~\ref{ass:common_ambiguity} hold.
If $F_f$ is monotone on $\mathcal X$, then
$\mathcal T_M$ is monotone on $\mathcal X\times\mathcal U
\times\mathbb R_{\geq0}^{2ML}.$
\end{proposition}
\emph{Proof:}
First, we note that each component of $c_M$ is convex in $(x,u)$. Specifically, 
$c_{k\ell}^{(1)}$ is a sum of a  convex function with affine terms, while $c_{k\ell}^{(2)}$ is the norm of an affine mapping
minus $\lambda$.
Since $\mu\geq0$, $L_M$ is convex in $(x,u)$
and affine in $\mu$.
Consequently, its saddle operator
\begin{equation}
(x,u,\mu)\mapsto
\begin{bmatrix}
\partial_{(x,u)}L_M(x,u,\mu)\\[1mm]
-c_M(x,u)
\end{bmatrix}.
\end{equation}
is monotone on $\mathcal X\times\mathcal U
\times\mathbb R_{\geq0}^{2ML}.$
Adding the monotone mapping $\operatorname{col}(F_f(x),0,0)$ preserves
monotonicity and concludes the proof. \hfill $\blacksquare$

\begin{remark}
If $r_M$ is convex and
$\prod_{i\in\mathcal N}\partial_{x_i}r_M(x)
=\partial r_M(x)$ on $\mathcal X$, then it holds that
$F_K=F_f+\partial r_M$ is monotone whenever $F_f$
is monotone. This condition, however, holds when $r_M$ is
differentiable or separable across agents' decisions,
but does not follow from homogeneous ambiguity alone. \hfill $\square$
\end{remark}
 
 This conclusion does not automatically
extend to the heterogeneous formulation with private
epigraphic constraints. The next subsection presents
a computational scheme for that general formulation.

\subsection{Iterative scheme for the finite-dimensional DRGNEP}
We consider a primal-dual scheme for computing equilibria of the
finite-dimensional reformulation in
Proposition~\ref{lem:polytopic}. In the general heterogeneous case,
the private epigraphic constraints do not directly yield the
standard monotonicity conditions required for a convergence
analysis. Hence, the method is used here as a computational scheme. 
For the polytopic reformulation used in this subsection, we denote 
\[
\gamma_i:=
\col\bigl((\gamma_{k_i,\ell_i})_{k_i,\ell_i}\bigr)
\in\mathbb R_{\geq0}^{mK_iL_i}, \ 
y_i:=\col(x_i,\lambda_i,\gamma_i,s_i),
\]
and $Y_i:=
X_i\times\mathbb R_{\geq0}
\times\mathbb R_{\geq0}^{mK_iL_i}
\times\mathbb R^{K_i}$,
where $m$ is the number of rows of $C_i$.
These definitions of $y_i$ and $Y_i$ are used
throughout Algorithm 1. To state it compactly, we denote:
$J_i(y_i)=\lambda_i \epsilon_i + \frac{1}{K_i}\sum_{k_i=1}^{K_i} s_{k_i}$, while
$$g^{(1)}_{k_i\ell_i} = b_{\ell_i}(x)+a_{\ell_i}^\top(x)\xi_i^{(k_i)}+ \gamma_{k_i, \ell_i}^\top (d_i-C_i\xi^{(k_i)}_i)- s_{k_i},$$
$$g^{(2)}_{k_i\ell_i}=\|C^\top_i\gamma_{k_i, \ell_i}-a_{\ell_i}(x)\| - \lambda_i,$$
and $g_i(y_i, x_{-i})=\col(\col(([ g^{(1)}_{k_i\ell_i}; \ g^{(2)}_{k_i\ell_i}]))_{k_i=1}^{K_i})_{\ell_i=1}^{L_i}$.
To include the shared constraints, define
\[
\widetilde g_i(y_i,x_{-i})
:=
\operatorname{col}\bigl(g_i(y_i,x_{-i}),\psi(x)\bigr).
\]
Let $n_{\mu_i}:=2L_iK_i+\tilde n$ and let
$\mu_i\in\mathbb R_{\geq0}^{n_{\mu_i}}$
be the associated multiplier vector. At iteration $k$, for each component
$\widetilde g_{i,r}$ of $\widetilde g_i$, select
\begin{equation}
v_{i,r}^{(k)}
\in
\partial_{y_i}\widetilde g_{i,r}
(\bar y_i^{(k)},x_{-i}^{(k)}), \nonumber 
\end{equation}
where $r=1,\ldots,n_{\mu_i}$. 
For this scheme, assume that
$a_{\ell_i}(\cdot,x_{-i})$ is affine and
$b_{\ell_i}(\cdot,x_{-i})$ is finite and convex
on an open convex neighborhood of $X_i$. For the Euclidean norm constraints, we can denote 
$u'=C_i^\top\gamma_{k_i,\ell_i}-a_{\ell_i}(x)$, and choose
$q=u'/\|u'\|$ if $u'\neq0$, and $q=0$ otherwise. The vector $u'$ is evaluated at
$(\bar y_i^{(k)},x_{-i}^{(k)})$.
 Here, $\delta\in[0,1)$ is the inertial parameter used to damp oscillations observed numerically. The physical constraints defining $\mathcal X$ are included in $\widetilde g_i$, whereas the Wasserstein epigraphic constraints remain private to each agent. Since Algorithm~1 uses individual multipliers, it is designed for a general DRGNE, and is therefore used as a heuristic. Requiring common multipliers for the physical shared
constraints imposes an additional equilibrium selection
condition.
 In the next section, we instead focus on variational DRGNE  and develop an equilibrium seeking method with convergence guarantees leveraging their structure. 

\section{Monotonicity, Sensitivity and vDRGNE Seeking}
This section establishes sufficient conditions under which
the robust pseudo-subgradient mapping remains strongly monotone under heterogeneous ambiguity. Furthermore, we derive
an equilibrium sensitivity bound and analyze the convergence of a projected
method with inexact robust gradient evaluations.
\subsection{Monotonicity under heterogeneous ambiguity}
\label{sec:heterogeneous_monotonicity}
Section IV-D establishes monotonicity of an extended
saddle mapping under a common ambiguity set and a shared
uncertain cost term. Here, we study the original mapping
$\mathcal F_K$ under heterogeneous ambiguity sets. Monotonicity of $F_K$ is not guaranteed in general under
heterogeneous ambiguity sets. The results below
identify additional structural and regularity
conditions under which monotonicity, or strong
monotonicity, is preserved. For costs of the form
$h_i(x,\xi_i)=f_i(x)+r_i(x_i,\xi_i)$, strong monotonicity is
immediately preserved. Indeed, if $F_f$ is $\alpha$-strongly monotone
and $r_i(\cdot,\xi_i)$ is $\mu_i$-strongly convex, then
$\mathcal F_K$ is $(\alpha+\min_i\mu_i)$-strongly monotone.
We next consider a more general case where uncertainty interacts with all agents'
decisions.
To this end, let $G(x,\xi)
:= \operatorname{col}\left(
(\nabla_{x_i}h_i(x,\xi_i))_{i\in\mathcal N}
\right)$,
where $\xi:=\operatorname{col}((\xi_i)_{i\in\mathcal N})$. We now impose the following assumption:

\begin{algorithm}[t] 
\caption{Primal-dual algorithm for the lifted \gls{DRGNEP}}
\label{Algorithm_apriori}
\begin{algorithmic}[1]
  \STATE \textbf{Initialization:} Choose $\delta\in[0,1)$ and set \( y^{(0)}_i \in Y_i, \  \xi^{(k_i)}_i \in \Xi^{K_i}_i, \mu_i^{(0)}\in\mathbb R_+^{n_{\mu_i}}\),
$\bar y_i^{(-1)}=y_i^{(0)}$,
$\bar\mu_i^{(-1)}=\mu_i^{(0)}$, $\epsilon_i \in \mathbb{R}_{\geq 0}, \tau_i > 0$ for all \( i \in \mathcal{N} \).
  \STATE \textbf{Iteration} \( k \)
  \STATE Agent \( i \):
  \STATE (1) Updates the variables:
  \STATE \quad \( \bar{y}^{(k)}_i = (1-\delta) y^{(k)}_i + \delta \bar{y}^{(k-1)}_i \)
  \STATE \quad \( \bar{\mu}^{(k)}_i = (1-\delta) \mu^{(k)}_i + \delta \bar{\mu}^{(k-1)}_i \)
  \STATE (2) Receives \( x^{(k)}_j \) for all \( j \neq i \), then updates:
\STATE \quad \( y^{(k+1)}_i = \operatorname{proj}_{Y_i} \big[ \bar{y}^{(k)}_i - \tau_i \big( \nabla_{y_i} J_i(y^{(k)}_i) + \sum\limits_{r=1}^{n_{\mu_i}}\mu_{i,r}^{(k)}v_{i,r}^{(k)} \big) \big] \) \label{update_step}
\STATE \quad \( \mu^{(k+1)}_i = \operatorname{proj}_{\mathbb{R}^{n_{\mu_i}}_{\geq 0}} \big[ \bar{\mu}^{(k)}_i + \tau_i \widetilde{g}_i(y^{(k)}_i, x^{(k)}_{-i}) \big] \)
  \STATE \textbf{Until convergence}
\end{algorithmic}
\end{algorithm}
\begin{assumption} \label{ass:monotonicity}
There exists a scalar $\alpha\geq0$ such that 
$(x-y)^\top
\bigl(
G(x,\xi)-G(y,\xi)
\bigr)
\geq
\alpha\|x-y\|^2$,
for all $x,y\in \mathcal{X}$ and
$\xi\in\prod_{i\in\mathcal N}\Xi_i$. \hfill $\square$
\end{assumption}
Lemma~\ref{lem:robust_subgradient} represents robust
partial subgradients through worst-case distributions.
We now impose a bound on the variation of distributions
representing these subgradients as the decision changes.
\begin{assumption}\label{assum:adversarial_sensitivity}
For every $i\in\mathcal N$, there exists
$\kappa_{i,K}\geq0$ such that, for all $x,y\in X$, with
$x\neq y$, any subgradients $v_i\in\partial_{x_i}\mathcal{R}_i(x)$ and
$w_i\in\partial_{x_i}\mathcal{R}_i(y)$, there exist probability measures $\mathbb Q_i^x\in\mathcal A_{i,K}(x),
\mathbb Q_i^y\in\mathcal A_{i,K}(y)$ for which it holds that $v_i=
\mathbb E_{\mathbb Q_i^x}
[\nabla_{x_i}h_i(x,\xi_i)],
w_i=
\mathbb E_{\mathbb Q_i^y}
[\nabla_{x_i}h_i(y,\xi_i)]$,
and
$d_W(\mathbb Q_i^x,\mathbb Q_i^y)
\leq
\kappa_{i,K}\|x-y\|$.
\hfill $\square$
\end{assumption}

Assumption~\ref{assum:adversarial_sensitivity} imposes a Lipschitz regularity condition,
in the Wasserstein metric, on active distributions
representing arbitrarily selected robust subgradients.
It holds with $\kappa_{i,K}=0$ whenever the subgradients
$v_i$, $w_i$ can be represented by a common
worst-case distribution, for every
$x,y\in\mathcal X$ and every
$v_i\in\partial_{x_i}\mathcal R_i(x)$,
$w_i\in\partial_{x_i}\mathcal R_i(y)$.
A simpler sufficient condition is that the worst-case
distribution denoted by the singleton $\{\mathbb Q_i^\star(x)\}
:=
\operatorname*{arg\,max}_{
\mathbb Q_i\in\mathbb B_{\epsilon_i}(\hat{\mathbb P}_{K_i})}
\mathbb E_{\mathbb Q_i}[h_i(x,\xi_i)]$, is unique for every $x\in\mathcal X$
and that the mapping
$x\mapsto\mathbb Q_i^\star(x)$ is Lipschitz continuous
in the Wasserstein metric, i.e., 
\[
d_W\!\left(\mathbb Q_i^\star(x),\mathbb Q_i^\star(y)\right)
\le\kappa_{i,K}\|x-y\|,
\text{ for all } x,y\in\mathcal X, 
\]
and for each $i\in\mathcal N$.
In general, such assumption can be restrictive when
worst-case distributions are nonunique or vary
discontinuously with the decision profile.

For ease of notation, let us now define the quantities $L_i:=
\big(
\sum_{j=1}^n L_{i,j}^2
\big)^{1/2}$, and
$\eta_K:=
\Big(
\sum_{i=1}^N
L_i^2\kappa_{i,K}^2
\Big)^{1/2}$. Then, the following result holds.

\begin{theorem}\label{thm:general_monotonicity}
Let Assumptions~1-4,
\ref{assum:compact_uncertainty},
\ref{ass:monotonicity}, and
\ref{assum:adversarial_sensitivity} hold. 
Then, for every $x,y\in \mathcal{X}$,
$v\in\mathcal F_K(x)$ and
$w\in\mathcal F_K(y)$,
\[
(x-y)^\top(v-w)
\geq
(\alpha-\eta_K)\|x-y\|^2.
\]
Consequently,  $F_K$
is monotone if $\alpha\geq\eta_K$ and strongly
monotone with modulus $\alpha-\eta_K$ if
$\alpha>\eta_K$. \hfill $\blacksquare$
\end{theorem}
\begin{proof}
 Since the result is immediate if $x=y$, we focus on the case where $x\neq y$. For $x\neq y$, Lemma~\ref{lem:robust_subgradient}
and Assumption~\ref{assum:adversarial_sensitivity}
provide worst-case distributions
$\mathbb Q_i^x\in\mathcal A_{i,K}(x)$ and
$\mathbb Q_i^y\in\mathcal A_{i,K}(y)$ representing the selected robust
subgradients $v_i\in\partial_{x_i}\mathcal R_i(x)$ and
$w_i\in\partial_{x_i}\mathcal R_i(y)$, respectively, as characterized in
Lemma~\ref{lem:robust_subgradient}, and satisfying
\[
d_W(\mathbb Q_i^x,\mathbb Q_i^y)
\leq\kappa_{i,K}\|x-y\|.
\]
For ease of notation, denote $\mathbb Q^x:=\bigotimes_{i=1}^N\mathbb Q_i^x$.

Adding and subtracting
$\mathbb E_{\mathbb Q_i^x}[\nabla_{x_i}h_i(y,\xi_i)]$
in each block yields:
\begin{align*}
&(x-y)^\top(v-w)
= \\
&\sum_{i=1}^N
(x_i-y_i)^\top
\mathbb E_{\mathbb Q_i^x}
\!\left[
\nabla_{x_i}h_i(x,\xi_i)
-
\nabla_{x_i}h_i(y,\xi_i)
\right]
\\
&+
\sum_{i=1}^N
(x_i-y_i)^\top
\left(
\mathbb E_{\mathbb Q_i^x}
[\nabla_{x_i}h_i(y,\xi_i)]
-
\mathbb E_{\mathbb Q_i^y}
[\nabla_{x_i}h_i(y,\xi_i)]
\right)
\\
&\overset{(a)}{\geq}
\alpha\|x-y\|^2
\\ 
&-\sum_{i=1}^N
\|x_i-y_i\|
\left\|
\mathbb E_{\mathbb Q_i^x}
[\nabla_{x_i}h_i(y,\xi_i)]
-
\mathbb E_{\mathbb Q_i^y}
[\nabla_{x_i}h_i(y,\xi_i)]
\right\|
\\
&\overset{(b)}{\geq}
\alpha\|x-y\|^2
-
\sum_{i=1}^N
L_i\|x_i-y_i\|
d_W(\mathbb Q_i^x,\mathbb Q_i^y)
\\
&\overset{(c)}{\geq}
\alpha\|x-y\|^2
-
\|x-y\|
\sum_{i=1}^N
L_i\kappa_{i,K}\|x_i-y_i\|
\\
&\overset{(d)}{\geq}
(\alpha-\eta_K)\|x-y\|^2,
\end{align*}
where (a) follows by integrating Assumption~\ref{ass:monotonicity} with respect to
$\mathbb{Q}^x$ and applying the Cauchy-Schwarz inequality, (b) follows from
Assumption~3 and the Kantorovich-Rubinstein theorem,
(c) follows from Assumption~\ref{assum:adversarial_sensitivity},
and (d) follows from the Cauchy-Schwarz inequality, i.e.,
$\sum_{i=1}^N
L_i\kappa_{i,K}\|x_i-y_i\|
\leq
\eta_K
\big(
\sum_{i=1}^N\|x_i-y_i\|^2
\big)^{1/2}
=
\eta_K\|x-y\|.$
\end{proof}
Note that if $\alpha<\eta_K$, we cannot know from this result alone whether the mapping is monotone or not. Furthermore, Theorem~\ref{thm:general_monotonicity} shows that the variation of the decision-dependent worst-case distributions reduces the monotonicity modulus $\alpha$ by at most $\eta_K := \sqrt{\sum_{i=1}^{N} L_i^2 \kappa_{i,K}^2}$. 

\subsection{Finite-sample equilibrium sensitivity}
The pseudogradient sensitivity result of
Lemma~\ref{rho_upper_bound} can now be used at the
equilibrium level by exploiting the strong monotonicity
established in Theorem~\ref{thm:general_monotonicity}.
\begin{corollary}
\label{cor:eq_sensitivity}
Suppose that the conditions of Lemma~\ref{rho_upper_bound}
and Theorem~\ref{thm:general_monotonicity} hold with
$\mu_K:=\alpha-\eta_K>0$. Then, for every
$x_K^\star\in\operatorname{SOL}(\mathcal{X},\mathcal F_K)$ and
$x_{\mathbb P}^\star\in
\operatorname{SOL}(\mathcal{X},F_{\mathbb P})$, it holds that: $\|x_K^\star-x_{\mathbb P}^\star\|
\le
\dfrac{\sqrt{\rho_K}}{\mu_K}$. \hfill $\blacksquare$
\end{corollary}
\begin{proof} 
Since for $x_K^\star=x_{\mathbb P}^\star$, the result is
immediate, we focus on the case where
$x_K^\star \neq x_{\mathbb P}^\star$. Since
$x_K^\star\in\operatorname{SOL}(\mathcal{X},\mathcal F_K)$, there
exists $v_K^\star\in\mathcal F_K(x_K^\star)$ such that
\[
(y-x_K^\star)^\top v_K^\star\geq0,
\text{ for all } y\in \mathcal{X}.
\]
Let $w_K\in\mathcal F_K(x_{\mathbb P}^\star)$ be an
arbitrary subgradient. Then, by
Theorem~\ref{thm:general_monotonicity}, it holds that:
\begin{align*}
\mu_K\|x_K^\star-x_{\mathbb P}^\star\|^2
&\leq
(x_K^\star-x_{\mathbb P}^\star)^\top
(v_K^\star-w_K)\\
&\overset{(a)}{\leq}
-(x_K^\star-x_{\mathbb P}^\star)^\top w_K\\
&\overset{(b)}{\leq}
(x_K^\star-x_{\mathbb P}^\star)^\top
\bigl(F_{\mathbb P}(x_{\mathbb P}^\star)-w_K\bigr)\\
&\overset{(c)}{\leq}
\|x_K^\star-x_{\mathbb P}^\star\|
\|F_{\mathbb P}(x_{\mathbb P}^\star)-w_K\|.
\end{align*}
Here, (a) follows from the robust VI with
$y=x_{\mathbb P}^\star$, (b) follows from the nominal VI with $y=x_K^\star$, and (c) follows from
Cauchy-Schwarz inequality. Dividing by
$\mu_K\|x_K^\star-x_{\mathbb P}^\star\|>0$ and taking the
infimum over $w_K\in\mathcal F_K(x_{\mathbb P}^\star)$
yields
\[
\|x_K^\star-x_{\mathbb P}^\star\|
\leq
\frac{
\operatorname{dist}\!\left(
F_{\mathbb P}(x_{\mathbb P}^\star),
\mathcal F_K(x_{\mathbb P}^\star)
\right)}
{\mu_K}.
\]
Finally, evaluating the uniform pseudo-subgradient deviation bound of
Lemma~\ref{rho_upper_bound} at the nominal vSGNE
$x_\mathbb{P}^\star\in\mathcal X$ gives
\begin{align}
\operatorname{dist}
\bigl(
F_{\mathbb P}(x_\mathbb{P}^\star),
\mathcal F_K(x_\mathbb{P}^\star)
\bigr)
&\leq
\operatorname{dev}_K(x_\mathbb{P}^\star) \nonumber 
\\
&\leq
\|\mathcal F_K-F_{\mathbb{P}}\|_{\mathcal X}
\leq\sqrt{\rho_K}, \nonumber 
\end{align}
which concludes the proof.
\end{proof}

\begin{remark}
 The bound in Corollary~2 depends on
$d_W(\hat{\mathbb{P}}_{K_i},\mathbb{P}_i)$ and therefore on the collection of unknown true
distributions of the agents. A confidence calibrated bound follows directly
from Lemma~4. Specifically, if it holds that
$d_W(\hat{\mathbb{P}}_{K_i},\mathbb{P}_i)\leq\epsilon_i(K_i,\beta_i)$ for all
$i\in\mathcal N$, it holds that
$\widetilde w_i\leq2\epsilon_i(K_i,\beta_i)$ and 
\[
\|x_K^\star-x_\mathbb{P}^\star\|
\leq \textstyle 
\frac{2}{\mu_K}
\sqrt{
\sum_{i=1}^N\sum_{j=1}^n
L_{i,j}^2\epsilon_i(K_i,\beta_i)^2
}.
\]
By Lemma~4 and Bonferroni's inequality, this bound then holds with
probability at least $1-\sum_{i\in\mathcal N}\beta_i$.
\end{remark}
\subsection{Inexact vDRGNE Seeking}
\label{sec:inexact_wcg}
Algorithm 1 operates on the explicit lifted GNEP
and targets general DRGNEs through the equilibrium
correspondence in Lemma 9. However,  no general
convergence guarantee is provided for that scheme and it therefore acts as a heuristic.
Algorithm 2, proposed in this section, instead operates directly on the original
robust VI and targets a vDRGNE.
Its convergence analysis requires  certain properties of the robust mapping such as single-valuedness, Lipschitz continuity and
strong monotonicity. This additional structure, however, allows us to establish convergence under vanishing oracle errors, an asymptotic error bound under bounded oracle errors, and linear convergence under exact oracle evaluations.
\GP{\begin{algorithm}[t]
\caption{Inexact Projected Gradient vDRGNE Seeking}
\label{alg:inexact_wcg}
\begin{algorithmic}[1]
\STATE \textbf{Initialization:}
Choose $x^{(0)}\in\mathcal X$, $\gamma>0$,  
$\xi_i^{K_i}\in\Xi_i^{K_i}$,
$\epsilon_i\in\mathbb R_{\geq0}$ for each $i\in\mathcal N$.
\STATE \textbf{Iteration $k$}
\STATE \textbf{Agent $i$:}
\STATE (1) Receives $x_j^{(k)}$, for all $j\neq i$ then updates:
\STATE (2) Computes an inexact gradient
$\hat g_i^{(k)}$ satisfying
\[
\|\hat g_i^{(k)}-g_i^{(k)}\|
\leq\delta_i^{(k)},
g_i^{(k)}
\in\partial_{x_i}\mathcal{R}_i(x^{(k)}).
\]
\STATE (3) Updates
\[
x^{(k+1)}
=
\operatorname{proj}_{\mathcal{X}}
[
x^{(k)}-\gamma\hat g^{(k)}].
\]
\STATE \textbf{Until convergence.}
\end{algorithmic}
\end{algorithm}}
The previous results characterize when the original robust
pseudogradient remains strongly monotone despite heterogeneous
ambiguity. We now exploit this property to devise an
equilibrium-seeking method that operates directly on the
DRVI, without introducing the auxiliary variables of the
finite-dimensional GNEP in Section~IV. By Lemma~\ref{lem:robust_subgradient}, for every decision  $x\in \mathcal{X}$ and
$i\in\mathcal N$, any active worst-case distribution
$\mathbb Q_i\in\mathcal A_{i,K}(x)$ generates a robust
subgradient
\begin{align}
g_i(x,\mathbb Q_i)
:=
\mathbb E_{\mathbb Q_i}
[\nabla_{x_i}h_i(x,\xi_i)]
\in\partial_{x_i}\mathcal{R}_i(x).
\end{align}
We allow this quantity to be computed inexactly, i.e, 
at iteration $k$, agent $i$ computes $\hat g_i^{(k)}$
such that, for some
$g_i^{(k)}\in\partial_{x_i}\mathcal{R}_i(x^{(k)})$, with  $\|\hat g_i^{(k)}-g_i^{(k)}\|
\leq\delta_i^{(k)}$,
where $\delta_i^{(k)}\geq0$ denotes an approximation
error. 
Since $\mathcal X$ contains shared coupling constraints, the projection
onto $\mathcal X$ is performed as a coordinated collective update and
is not generally separable across agents. Future work will focus on the development of a distributed version of Algorithm 2. At iteration $k$, suppose that the inner \gls{DRO} oracle returns a
primal feasible distribution
$\widetilde Q_i^{(k)}\in
\mathbb{B}_{\epsilon_i}(\hat{\mathbb{P}}_{K_i})$
and a certified upper bound
$\mathcal{U}_i^{(k)}\ge \mathcal{R}_i(x^{(k)})$.
Define $\mathcal{L}_i^{(k)}
:=
\mathbb E_{\widetilde Q_i^{(k)}}
[h_i(x^{(k)},\xi_i)]$,
and 
$\Delta_i^{(k)}
:=
\mathcal{U}_i^{(k)}-\mathcal{L}_i^{(k)}$.
Then, by primal and dual feasibility, it holds that:
\[
0\le
\mathcal{R}_i(x^{(k)})-\mathcal{L}_i^{(k)}
\le
\Delta_i^{(k)}.
\]
For the Wasserstein ambiguity sets considered here,
$\mathcal{U}_i^{(k)}$ can be obtained from any dual feasible point of
the reformulation in Theorem~\ref{reformulation_1}, while $\mathcal{L}_i^{(k)}$ is provided
by any primal feasible transport plan. Then, the following result connects the chosen accuracy to the error bound of the approximated subgradient.

\begin{lemma}
\label{lem:value_gap_oracle}
For a given agent $i\in\mathcal N$ and decision $x\in\mathcal X$. Suppose that
$\mathcal{R}_i(\cdot,x_{-i})$ has an $\ell_i$-Lipschitz gradient locally on a Euclidean ball 
$\mathbb{B}(x_i,r_i)$. Let
$\widetilde Q_i\in\mathbb B_{\epsilon_i}(\hat{\mathbb P}_{K_i})$
be an approximate worst-case distribution returned by the inner \gls{DRO} oracle and satisfy
\[
\mathcal{R}_i(x)-\mathbb E_{\widetilde Q_i}[h_i(x,\xi_i)]\le \Delta_i,
\]
and consider  $\widehat g_i
:=
\mathbb E_{\widetilde Q_i}
[\nabla_{x_i}h_i(x,\xi_i)]$.
Then, if $\Delta_i\le \dfrac{\ell_i r_i^2}{2}$,
it holds that $\|\widehat g_i-\nabla_{x_i}R_i(x)\|
\le
\sqrt{2\ell_i\Delta_i}.$ 
\end{lemma}
\emph{Proof:}
For any $\|d\|\le r_i$, feasibility of $\widetilde Q_i$ and convexity
in $x_i$ give $R_i(x_i+d,x_{-i})
\ge
R_i(x)-\Delta_i+\widehat g_i^\top d$.
By $\ell_i$-smoothness, it then holds that
\[
\mathcal{R}_i(x_i+d,x_{-i})
\le
\mathcal{R}_i(x)+\langle\nabla_{x_i}\mathcal{R}_i(x),d\rangle
+\dfrac{\ell_i}{2}\|d\|^2.
\]
Hence, the two inequalities imply that: 
\[
\langle\widehat g_i-\nabla_{x_i}\mathcal{R}_i(x),d\rangle
\le
\Delta_i+\dfrac{\ell_i}{2}\|d\|^2.
\]
Choosing $d=t e/\|e\|$, where
$e=\widehat g_i-\nabla_{x_i}\mathcal{R}_i(x)$, yields $\|e\|
\le
\dfrac{\Delta_i}{t}+\dfrac{\ell_i t}{2}, 0<t\le r_i.$
Taking the infimum and setting
$t=\sqrt{2\Delta_i/\ell_i}$ concludes the proof. \hfill $\blacksquare$

For the next result, we denote, for ease of notation,  $\hat g^{(k)}
:=
\operatorname{col}
((\hat g_i^{(k)})_{i\in\mathcal N})$, and  
$\delta^{(k)}
:=
\big(
\sum_{i\in\mathcal N}(\delta_i^{(k)})^2
\big)^{1/2}.$

\begin{theorem}
\label{thm:inexact_wcg}
Suppose that the conditions of
Theorem~\ref{thm:general_monotonicity} hold with $\mu_K:=\alpha-\eta_K>0.$
Moreover, suppose that $\mathcal F_K$ is single-valued, denoted by $F_K$, and $L_K$-Lipschitz continuous on $\mathcal{X}$.
Let $0<\gamma<\frac{2\mu_K}{L_K^2}$
and define $q_K:=
\sqrt{1-2\gamma\mu_K+\gamma^2L_K^2}<1$.
Then, the sequence generated by
Algorithm~\ref{alg:inexact_wcg} satisfies 
\begin{align}
\|x^{(k+1)}-x_K^\star\|
\leq
q_K\|x^{(k)}-x_K^\star\|
+\gamma\delta^{(k)},
\end{align}
where $x_K^\star$ is the unique vDRGNE, i.e., the unique solution of
$\operatorname{VI}(\mathcal{X},F_K)$. 
\end{theorem}
\begin{proof}
Since $\mu_K>0$, Theorem~\ref{thm:general_monotonicity}
implies that $F_K$ is strongly monotone and therefore the
solution $x_K^\star$ is unique. Moreover, since
$x_K^\star$ solves $\operatorname{VI}(\mathcal{X},F_K)$, it holds that
$x_K^\star
=
\operatorname{proj}_{\mathcal{X}}
\left[
x_K^\star-\gamma F_K(x_K^\star)
\right].$
Let $e^{(k)}
:=
\hat g^{(k)}-F_K(x^{(k)})$.
Since $F_K$ is single-valued, the oracle condition yields
$\|e^{(k)}\|\leq\delta^{(k)}$. Therefore, we have the following inequalities:
\begin{align*}
&\|x^{(k+1)}-x_K^\star\| \\
&\overset{(a)}{\leq}
\bigl\|
x^{(k)}-x_K^\star
-\gamma
\bigl(F_K(x^{(k)})-F_K(x_K^\star)\bigr)
-\gamma e^{(k)}
\bigr\|
\\
&\overset{(b)}{\leq}
\bigl\|
x^{(k)}-x_K^\star
-\gamma
\bigl(F_K(x^{(k)})-F_K(x_K^\star)\bigr)
\bigr\|
+\gamma\delta^{(k)}
\\
&\overset{(c)}{\leq}
q_K\|x^{(k)}-x_K^\star\|
+\gamma\delta^{(k)}.
\end{align*}
Here, (a) follows from the nonexpansiveness of the
projection, (b) follows from the triangle inequality and the
oracle error bound, while (c) follows from
$\mu_K$-strong monotonicity and $L_K$-Lipschitz continuity
of $F_K$, since
\begin{align*}
&\bigl\|
x^{(k)}-x_K^\star
-\gamma
(F_K(x^{(k)})-F_K(x_K^\star))
\bigr\|^2
\\
&\qquad\leq
\left(
1-2\gamma\mu_K+\gamma^2L_K^2
\right)
\|x^{(k)}-x_K^\star\|^2, 
\end{align*}
thus concluding the proof. 
\end{proof}
Consequently, if $\delta^{(k)}\to 0$, then
$x^{(k)}\to x_K^\star$. If instead
$\delta^{(k)}\leq\bar\delta$ for all $k$, then
\[
\limsup_{k\to\infty}\|x^{(k)}-x_K^\star\|
\leq
\frac{\gamma\bar\delta}{1-q_K}.
\]
In particular, for exact oracle evaluations,
$\delta^{(k)}=0$, Algorithm~2 converges linearly to
$x_K^\star$ with contraction factor $q_K$.
\section{Numerical simulations}
We illustrate the finite-dimensional reformulation and Algorithm~1 on a peer-to-peer electricity market, and the monotonicity results and Algorithm~2 on a stochastic Nash-Cournot game. In the numerical experiments, the Wasserstein radii are prescribed
exogenously rather than calibrated through~(8), since our objective is
to illustrate the effect of heterogeneous ambiguity and varying levels
of distributional robustness on the equilibrium and the proposed
algorithms.

\subsection{Peer-to-peer electricity markets} \label{P2P_section}
We consider a \gls{DR} extension of the peer-to-peer
electricity market in~\cite{LeCadre2020}. Each participant
$i\in\mathcal N$ chooses generation $G_i$, demand $D_i$,
and bilateral trades
$\nu_i:=\{\nu_{mi}\}_{m\in\Omega_i}$, subject to local
bounds, trading capacities, reciprocity constraints, and
power balance, where $\Delta G_i$ denotes renewable
generation. As in~\cite{LeCadre2020}, reciprocity is written
as $\nu_{mi}\leq-\nu_{im}$. Since
uncertainty affects the generation cost parameters, each
participant constructs an individual Wasserstein ambiguity
set from private samples and radius $\epsilon_i$ and solves:

\begin{align*}
\begin{cases}
\begin{aligned}
\min_{G_i,D_i,\nu_i}\ \max_{Q_i \in \mathcal{P}_i}
\quad &
\mathbb{E}_{Q_i}[
J_{g_i}(G_i,\xi_i)
+ J_{p2p,i}(\nu_i)
+ J_{d_i}(D_i)]
\\
\text{s.t.}\quad
& \underline{D}_i \le D_i \le \overline{D}_i, \\
& \underline{G}_i \le G_i \le \overline{G}_i, \\
& \nu_{mi} \le \chi_{mi},
\quad \forall m \in \Omega_i,\ m \neq i, \\
& \nu_{mi} \le -\nu_{im},
\quad \forall m \in \Omega_i,\ m \neq i, \\
& D_i = G_i + \Delta G_i + \textstyle \sum_{m \in \Omega_i} \nu_{mi}.
\end{aligned}
\end{cases}
\end{align*}
\begin{figure}[t]
    \centering
    \includegraphics[width=0.75\columnwidth]
    {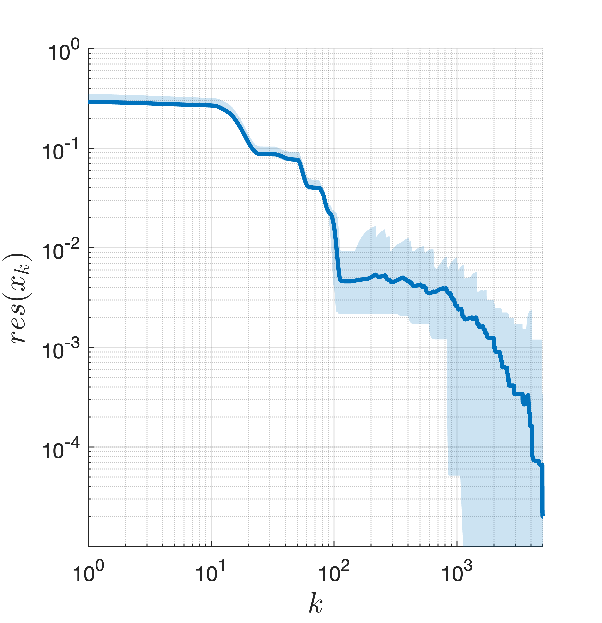}
    \caption{Empirical residual of Algorithm~1
    for $30$ P2P instances with heterogeneous multisamples,
    Wasserstein radii, and renewable generation values
    $\Delta G_i$. The solid line is the mean across
    simulations and the shaded region the observed range.}
    \label{P2P}
\end{figure}
We consider the generation cost
$J_{g_i}(G_i,\xi_i)=\xi_{i,1}G_i+\xi_{i,2}$,
where $\xi_i=(\xi_{i,1},\xi_{i,2})$ is a private uncertain parameter following an unknown probability distribution $\mathbb{P}_i$. The peer-to-peer trading cost is
$J_{p2p,i}(\nu_i)=\sum_{m \in \Omega_i} c_{mi}\nu_{mi}$,
where $c_{mi}$ denotes the bilateral trading price. The demand cost is
$J_{d_i}(D_i)=\omega_{i,1}(D_i-D_i^\ast)^2-\omega_{i,2}$,
where $D_i^\ast$ is the preferred demand and $\omega_{i,1},\omega_{i,2}$ are given parameters. Each participant constructs the individual ambiguity set
$\mathcal{P}_i=
\mathbb{B}_{\epsilon_i}(\hat{\mathbb{P}}_{K_i})$
from its private samples. We leverage Proposition~\ref{lem:polytopic} to obtain the
finite-dimensional reformulation and apply Algorithm~1.
We consider $N=5$ participants and 30 instances
with different private multisamples, heterogeneous
Wasserstein radii $\epsilon_i\in[0.1,0.25]$, and perturbed
renewable generation values $\Delta G_i$. For this
two-dimensional uncertainty, the radii are chosen as
robustness parameters independently of the statistical
calibration in \eqref{epsilon_i}, to illustrate
algorithmic behavior under heterogeneous ambiguity. Figure~\ref{P2P} shows  the residual iterates across the simulated
instances.

\subsection{Heterogeneous Nash--Cournot games} \label{NC_section}

We consider $N$ firms selecting quantities $x_i\in X_i$ under the shared market constraint $v(x)\leq0$.  Each firm solves
\begin{equation}
\left\{
\begin{aligned}
\min_{x_i\in X_i}\quad&
\sup_{\mathbb Q_i\in
\mathbb B_{\epsilon_i}(\hat{\mathbb P}_{K_i})}
\mathbb E_{\mathbb Q_i}
[J_i(x,\xi_i)]
\\
\text{s.t.}\quad&
v(x)\leq0 .
\end{aligned}
\right.
\label{eq:cournot}
\end{equation}
where
$J_i(x,\xi_i)=c_ix_i^2-(w_1-w_2\sum_{j\in\mathcal N}x_j)\xi_i x_i$.
The sample sizes $K_i$ and radii $\epsilon_i$ may differ across firms.
We next consider an instance of
\eqref{eq:cournot} to illustrate the results of
Section~\ref{sec:heterogeneous_monotonicity}. We consider $N=20$ firms with
$X_i=[\underline x_i,\overline x_i]$ and scalar compact
supports $\Xi_i=[\underline\xi_i,\overline\xi_i]$.
The firms are coupled through the aggregate capacity
constraint $\sum_{i=1}^{N}x_i\leq C.$ We choose $0<\underline x_i\le x_i\le\overline x_i$
and parameters satisfying
$w_1>w_2\sum_{j\in\mathcal N}\overline x_j$.
Hence, it holds that $x_i(w_2\sum_j x_j-w_1)<0$ throughout
the feasible set, so the worst-case distribution of
firm $i$ always minimizes the attainable mean of
$\xi_i$ and is independent of the decision profile. Hence, the worst-case distribution
of firm $i$ is characterized by the minimum attainable mean, i.e.,  $m_{i,K}
:=
\inf_{\mathbb Q_i\in
\mathbb B_{\epsilon_i}(\hat{\mathbb P}_{K_i})}
\mathbb E_{\mathbb Q_i}[\xi_i].$
For $\Xi_i=[\underline\xi_i,\overline\xi_i]$, this implies
$m_{i,K}
=
\max\{
\underline\xi_i,\,
\frac{1}{K_i}\sum_{k_i=1}^{K_i}\xi_i^{(k_i)}
-\epsilon_i
\}.$
Consequently, the robust pseudogradient is affine and takes the form:
\begin{align}
F_{K,i}(x)
=2c_i x_i-w_1m_{i,K}+
w_2m_{i,K}
\Big(
2x_i+{\textstyle \sum_{j\neq i}}x_j
\Big), \nonumber
\end{align}
or equivalently $F_K(x)=J_Kx+q_K$,
where $[J_K]_{ii}=2c_i+2w_2m_{i,K}$, and $[J_K]_{ij}=w_2m_{i,K}$, where $j\neq i$.
If all firms have the same worst-case mean,
$m_{i,K}=m_K$ and the cross derivatives are symmetric.
Under the natural conditions $c_i>0$, $w_2\ge0$,
and $m_K\ge0$, the resulting Jacobian is positive
definite, thus recovering the standard strongly
monotone Nash-Cournot structure.
Under heterogeneous ambiguity,
$[J_K]_{ij}=w_2m_{i,K}$ and
$[J_K]_{ji}=w_2m_{j,K}$ generally differ, so monotonicity is
no longer automatic. For a fixed realization $\xi$, the symmetric part of the Jacobian is $H_{ii}(\xi)=2c_i+2w_2\xi_i,
H_{ij}(\xi)=\dfrac{w_2}{2}(\xi_i+\xi_j)$.
 Consider the optimization problem of finding the smallest eigenvalue of $H$ as follows:
\[
\underline\alpha
:=
\min_{\xi\in\prod_{i\in\mathcal N}\Xi_i}
\lambda_{\min}(H(\xi)).
\]
If $\underline\alpha\ge0$, Assumption~\ref{ass:monotonicity} holds with
$\alpha=\underline\alpha$. Note that if $\underline\alpha>0$,
the corresponding mapping $G(\cdot,\xi)$ is strongly monotone
on $\mathcal X$ with monotonicity constant
$\underline\alpha$.
Since $x_i(w_2\sum_{j\in\mathcal N}x_j-w_1)<0$
for every $x\in\mathcal X$, the
 set of active distributions consists of those
attaining the minimum mean $m_{i,K}$ and does not
depend on $x$.
Moreover, the gradient depends on $\mathbb Q_i$
only through the mean $\mathbb E_{\mathbb Q_i}[\xi_i]$.
Hence, every active distribution produces the same
robust gradient, and a common active distribution
can be chosen at any two decision profiles.
Therefore, Assumption~\ref{assum:adversarial_sensitivity} holds with $\kappa_{i,K}=0$
for every $i\in\mathcal N$, and it holds that $\eta_K=0$.
For $\underline\alpha\ge0$,
Theorem~\ref{thm:general_monotonicity} therefore yields
$\mu_K=\underline\alpha$.
We apply Algorithm~2 with aggregate oracle
error $\delta^{(k)}=10^{-4}$ and consider
$\mu_K\in\{1,0.2,0.02,-0.4\}$.
For these experiments, we retain the notation
$\mu_K=\underline\alpha$. Note that the case $\mu_K=\underline\alpha=-0.4$ is included only as a numerical case study 
illustrating residual convergence even outside the assumptions
of Theorems \ref{thm:general_monotonicity} and \ref{thm:inexact_wcg}.

\section{Conclusion}
We studied data-driven Wasserstein \gls{DR} generalized Nash games with deterministic shared convex constraints and heterogeneous ambiguity sets built from private data. We established finite-sample guarantees and finite dimensional
reformulations for DRGNEs, while deriving sensitivity
bounds, almost-sure consistency, and, under certain conditions, convergence guarantees for vDRGNE seeking. Future research will develop scalable distributed primal-dual schemes that jointly handle physical coupling constraints and Wasserstein dual variables. Further directions include convergence under weaker monotonicity conditions, nonunique worst-case distributions as well as uncertain coupling constraints.

\begin{figure}[t]
    \centering
    \includegraphics[width=0.95\columnwidth]{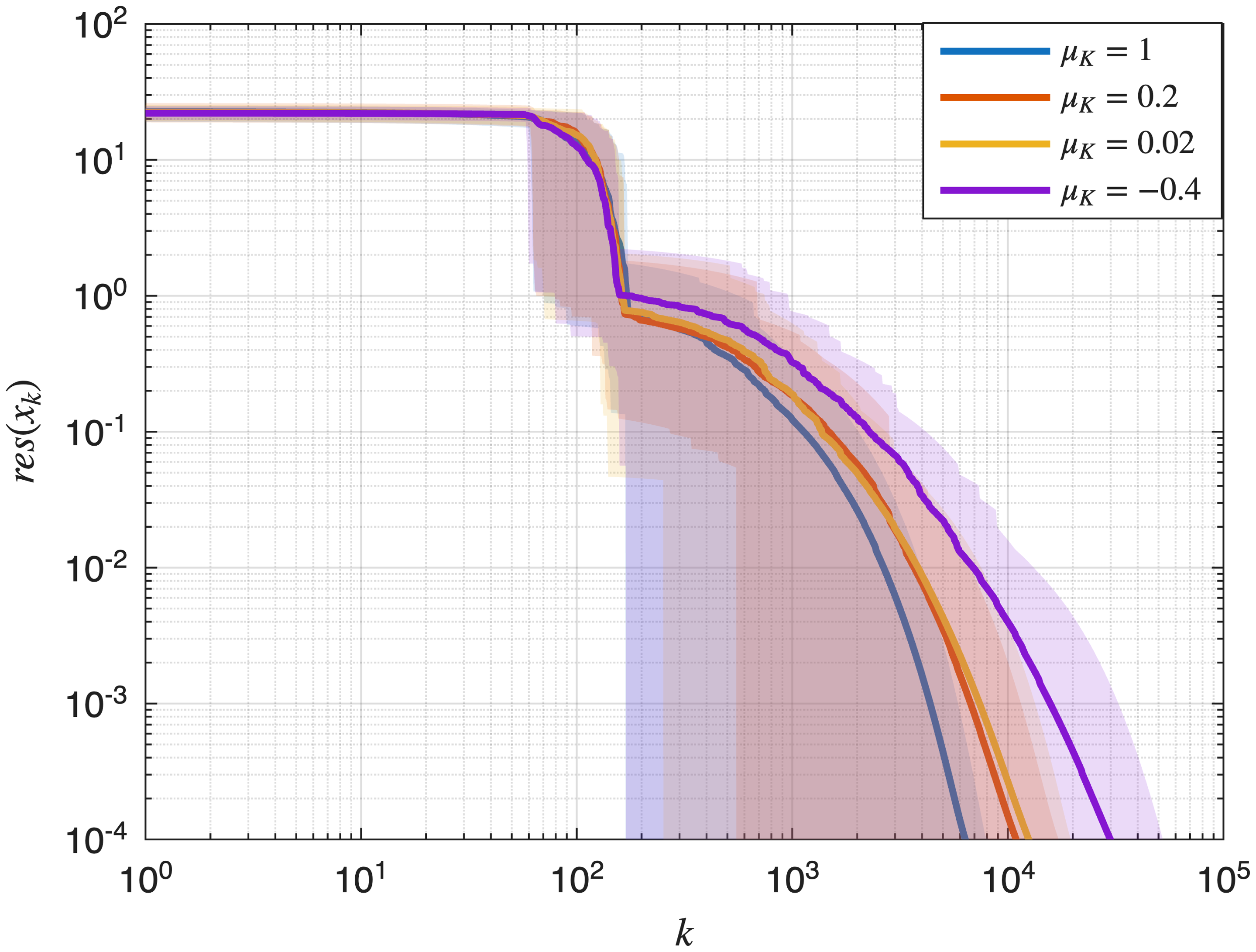}
    \caption{Residual of the inexact projected gradient method for the heterogeneous
    Nash-Cournot game. Solid lines are the mean over $50$
    randomized data and parameter realizations, while the shaded
    regions denote the observed range.}
    \label{fig:cournot_residual_inexact}
\end{figure}

\appendix

\subsection{Proofs of Section \ref{sec_problem}}\label{app_problem}
 
\begin{proof}[Proof of Lemma \ref{DRNE-DRVI}] Let
$x_K^\star\in
\operatorname{SOL}(\mathcal X,\mathcal F_K)$. By definition, there exists
$v_K^\star
=
\operatorname{col}
\left(
(v_{i,K}^\star)_{i\in\mathcal N}
\right)
\in\mathcal F_K(x_K^\star)$
such that:
\begin{align}
(y-x_K^\star)^\top v_K^\star\geq 0, \ \forall y\in\mathcal X.
\label{eq}
\end{align}
For a given $i\in\mathcal N$, consider any
$z_i\in\mathcal X_i(x_{-i,K}^\star)$ and the decision vector $z:=(z_i,x_{-i,K}^\star)\in\mathcal{X}$.
Substituting this point into
\eqref{eq} gives $(z_i-x_{i,K}^\star)^\top v_{i,K}^\star\geq 0$.
Moreover, since $v_K^\star\in\mathcal F_K(x_K^\star)$, we have that $v_{i,K}^\star
\in
\partial_{x_i}\mathcal R_i(x_K^\star)$.
By Assumption \ref{convex}, $\mathcal R_i(\cdot,x_{-i,K}^\star)$
is convex and the subgradient inequality gives:
\begin{align}
\mathcal R_i(z_i,x_{-i,K}^\star)
&\geq
\mathcal R_i(x_{i,K}^\star,x_{-i,K}^\star)
+
(z_i-x_{i,K}^\star)^\top v_{i,K}^\star
\nonumber  \\
&\geq
\mathcal R_i(x_{i,K}^\star,x_{-i,K}^\star). \nonumber
\end{align}
This implies that for each $i \in \mathcal{N}$ it holds that:
\begin{align}
x_{i,K}^\star
\in
\operatorname{argmin}\limits_{
x_i\in\mathcal X_i(x_{-i,K}^\star)}
\mathcal R_i(x_i,x_{-i,K}^\star), \nonumber 
\end{align}
i.e.,  $x_K^\star$ is a
\gls{DRGNE} of~\eqref{eq:dr}. Conversely, suppose that $\mathcal X=X$ and
$x_K^\star$ is thus a DRNE. For every $i\in\mathcal N$, the optimality condition means that there exists
$v_{i,K}^\star\in
\partial_{x_i}\mathcal R_i(x_K^\star)$ such that $(y_i-x_{i,K}^\star)^\top v_{i,K}^\star
    \geq 0$, for all  $y_i\in X_i$.
Stacking these subgradients and summing the individual
inequalities yields $(y-x_K^\star)^\top v_K^\star\geq 0, \forall y\in X$
for some $v_K^\star\in\mathcal F_K(x_K^\star)$, or in other words that
$x_K^\star\in\operatorname{SOL}(X,\mathcal F_K)$. 
\end{proof}

\subsection{Proofs of Section \ref{sec_prop}}\label{app_prop}
\begin{proof}[Proof of Lemma \ref{Lem:certificates_different}] 
Let $\xi^K:=(\xi_{K_i})_{i\in\mathcal N}
\in\prod_{i\in\mathcal N}\Xi_i^{K_i}$
denote the joint collection of multi-samples,
with probability distribution $\mathbb P^K$.
For each $i\in\mathcal N$, define the event
\[
C_i:=
\Big\{
\xi^K\in\prod_{j\in\mathcal N}\Xi_j^{K_j}:
d_W\!\left(
\mathbb P_i,\hat{\mathbb P}_{K_i}(\xi_{K_i})
\right)
\leq\epsilon_i(K_i,\beta_i)
\Big\}.
\]
If a multi-sample belongs to $C_i$, it holds that
$\mathbb P_i\in
\mathbb B_{\epsilon_i}(\hat{\mathbb P}_{K_i})$ and hence, by
the definition of $\mathcal R_i$, $\mathbb E_{\mathbb P_i}[h_i(x_K^*,\xi_i)]
\leq
\mathcal R_i(x_K^*;\hat{\mathbb P}_{K_i},
\epsilon_i(K_i,\beta_i))$.
Therefore, by Bonferroni's inequality and
Lemma~\ref{guarantees}:
\begin{align*}
\textstyle  \mathbb P^K\!\left(\bigcap_{i=1}^N C_i\right)
&\geq \textstyle 
1-\sum_{i=1}^N
\mathbb P^K(C_i^c)=
1-\sum_{i=1}^N
\mathbb P_i^{K_i}(C_i^c)\\
&\geq \textstyle 
1-\sum_{i=1}^N\beta_i,
\end{align*}
which concludes the proof.
\end{proof}
\vspace{-.5cm}
\begin{proof}[Proof of Lemma \ref{rho_convergence}] 
Assumption~\ref{light_tail} implies that $\mathbb P_i$ has a finite first moment. Hence, for every $i\in\mathcal N$, the empirical
distribution satisfies $d_W(\hat{\mathbb P}_{K_i},\mathbb P_i)
\rightarrow 0, 
\mathbb P^\infty\text{-almost surely}.$
Therefore, $\tilde w_i
= \epsilon_i(K_i,\beta_i^{(K_i)})
+
d_W(\hat{\mathbb P}_{K_i},\mathbb P_i)
\rightarrow 0$
almost surely for every $i\in\mathcal N$. Since $\mathcal N$ is
finite, these convergences hold simultaneously with
$\mathbb P^\infty$-probability one. Thus, $\rho_K
=
\sum_{i=1}^N\sum_{j=1}^n
L_{i,j}^2\tilde w_i^2
\rightarrow0$ a.s. as $\min_{i\in\mathcal N}K_i\to\infty$. 
\end{proof}

\subsection{Proofs of Section \ref{sec_finite}}\label{app_finite}
\begin{proof}[Proof of Theorem \ref{reformulation_1}]
For every $i\in\mathcal N$, for the $i$-th  problem in
\eqref{eq:dr},
the inner supremum is a Wasserstein \gls{DRO}
problem for agent $i$, where the strategic coupling with the other
players is retained through the dependence on $x_{-i}$ in
$h_i(x_i,x_{-i},\xi_i)$. By adapting \cite[Theorem 4.2]{mohajerin_esfahani_data-driven_2018}, and under
Assumption~\ref{pointwise_max}, it follows that
\begin{align}
&\sup_{\mathbb Q_i\in \mathbb B_{\epsilon_i}(\hat{\mathbb P}_{K_i})}
\mathbb E_{\mathbb Q_i}[h_i(x_i,x_{-i},\xi_i)]
= \nonumber \\
&\inf_{\lambda_i\ge 0}
\left\{
\lambda_i\epsilon_i
+
\frac{1}{K_i}\sum_{k_i=1}^{K_i}
\sup_{\xi_i\in\Xi_i}
(h_i(x_i,x_{-i},\xi_i)-\lambda_i\|\xi_i-\xi_i^{(k_i)}\|)
\right\}. \nonumber 
\end{align}
Under the dual attainment assumption, stated in the
theorem, player $i$'s problem is equivalently written as
\begin{align}
\min_{x_i\in \mathcal{X}_i(x_{-i}),\ \lambda_i\ge 0}
\left\{
\lambda_i\epsilon_i
+ \textstyle 
\frac{1}{K_i}\sum_{k_i=1}^{K_i}
g_i(x_i,x_{-i},\lambda_i,\xi_i^{(k_i)})
\right\}, \nonumber 
\end{align}
where
$g_i(x_i,x_{-i},\lambda_i,\xi_i^{(k_i)})
:=
\sup_{\xi_i\in\Xi_i}
(h_i(x_i,x_{-i},\xi_i)-\lambda_i\|\xi_i-\xi_i^{(k_i)}\|)$.
Since this reformulation holds for every agent $i\in\mathcal N$, the
game in \eqref{eq:dr} admits the stated equivalent reformulation. For $\epsilon_i>0$, the interchange above is justified by
\cite[Proposition 3.4]{shapiro_duality_2001}. For $\epsilon_i=0$, the Wasserstein ball reduces
to the singleton $\{\hat{\mathbb P}_{K_i}\}$, i.e., 
\[
\mathcal R_i(x;\hat{\mathbb P}_{K_i},0)
= \textstyle 
\frac{1}{K_i}\sum_{k_i=1}^{K_i}
h_i(x,\xi_i^{(k_i)}).
\]
Thus, for $\epsilon_i=0$, the empirical average can be used directly, without
introducing $\lambda_i$. 
\end{proof}
\begin{proof}[Proof of Lemma \ref{reformulation_2}]
Since
$-h_{\ell_i}(x_i,x_{-i},\cdot)+\iota_{\Xi_i}$
is proper, lower semicontinuous, and convex, while the norm
function is finite and continuous, Fenchel-Rockafellar duality yields $\sup_{\xi_i\in\Xi_i}
\min_{\|z_{k_i,\ell_i}\|_*\leq\lambda_i}
\phi_{k_i,\ell_i}(x_i,x_{-i},\xi_i,z_{k_i,\ell_i})
=
\min_{\|z_{k_i,\ell_i}\|_*\leq\lambda_i}
\sup_{\xi_i\in\Xi_i}
\phi_{k_i,\ell_i}(x_i,x_{-i},\xi_i,z_{k_i,\ell_i})$. Consequently, the original constraint is equivalent to the existence of
$z_{k_i,\ell_i}$ such that 
$\sup_{\xi_i\in\Xi_i}
\phi_{k_i,\ell_i}(x_i,x_{-i},\xi_i,z_{k_i,\ell_i})\leq s_{k_i}$, and
$\|z_{k_i,\ell_i}\|_*\leq\lambda_i$.
\end{proof}
\begin{proof}[Proof of Proposition \ref{lem:polytopic}]
For any arbitrary but fixed $x \in \mathcal{X}$, the finite maximum of affine
functions $h_i(x,\cdot)$ is Lipschitz continuous with constant
$L_i(x):=\max_{\ell_i}\|a_{\ell_i}(x)\|_*$.
For $\lambda_i\geq L_i(x)$, each inner supremum
in Theorem \ref{reformulation_1} equals the corresponding cost evaluated at the sample, since for each $k_i=1,\ldots,K_i$ it holds that: \[
\sup_{\xi_i\in\Xi_i}
\left\{
h_i(x,\xi_i)
-\lambda_i\|\xi_i-\xi_i^{(k_i)}\|
\right\}
=
h_i(x,\xi_i^{(k_i)}) .
\]
Thus, when $\epsilon_i=0$, the dual infimum is
attained. 
For $\epsilon_i>0$, choosing
$\xi_i=\xi_i^{(k_i)}$ in the supremum defining $g_i$
gives
\begin{multline*}
\lambda_i\epsilon_i+ \textstyle
\frac{1}{K_i}\sum_{k_i=1}^{K_i}
g_i(x_i,x_{-i},\lambda_i,\xi_i^{(k_i)}) \\
\geq
\lambda_i\epsilon_i+ \textstyle
\frac{1}{K_i}\sum_{k_i=1}^{K_i}
h_i(x,\xi_i^{(k_i)}).
\end{multline*}
Since the sample average is finite and independent
of $\lambda_i$, the right-hand side tends to
$+\infty$ as $\lambda_i\to+\infty$. Thus, the
dual objective is coercive on $[0,\infty)$.
Moreover, each $g_i$ is lower semicontinuous in
$\lambda_i$, being a supremum of continuous affine
functions of $\lambda_i$. The dual objective is
finite for $\lambda_i\geq L_i(x)$ and bounded
below by the sample average. Hence, it is proper
and lower semicontinuous. For $\epsilon_i>0$, it is coercive
in $\lambda_i$ and therefore attains its minimum on
$[0,\infty)$. The case $\epsilon_i=0$ is treated separately above.
Hence, the dual-attainment condition
required by Lemma \ref{reformulation_2} holds.
The player problems in \eqref{eq:dr} are coupled through both
the costs and the shared feasibility condition
$x_i\in\mathcal{X}_i(x_{-i})$. For fixed $i\in\mathcal N$ and
$x_{-i}$, convexity in $x_i$ and
\cite[Theorem 4.2, Corollary 5.1]
{mohajerin_esfahani_data-driven_2018} yield the reformulation
in \eqref{eq_prop}, while the feasibility condition remains
unchanged. Therefore, by the same argument as in
Lemma~\ref{equilibrium_equivalence}, any GNE
$(x^\ast_K,\lambda^\ast_K,\gamma^\ast_K,s^\ast_K)$ of
\eqref{eq_prop} has an $x$-component $x^\ast_K$ that is a
\gls{DRGNE} of \eqref{eq:dr}.
\end{proof}

\bibliographystyle{IEEEtran}
\bibliography{biblio_coop_barbara_v8}

\end{document}